\documentclass[11pt]{article}
\usepackage{}
\usepackage{mathrsfs}
\usepackage{amsfonts}
\usepackage{amsfonts,amssymb,mathrsfs}
\usepackage{amsmath,amscd}
\usepackage[dvips]{graphicx}
\usepackage[all]{xy}
\usepackage{graphicx}
\usepackage{fancyhdr}
\usepackage{mathptm,pslatex}
\usepackage{amsthm}

\begin{document}

\sloppy
\newtheorem{Def}{Definition}[section]
\newtheorem{Bsp}{Example}[section]
\newtheorem{Prop}[Def]{Proposition}
\newtheorem{Theo}[Def]{Theorem}
\newtheorem{Lem}[Def]{Lemma}
\newtheorem{Koro}[Def]{Corollary}
\theoremstyle{definition}
\newtheorem{Rem}[Def]{Remark}

\newcommand{\add}{{\rm add}}
\newcommand{\gd}{{\rm gl.dim }}
\newcommand{\dm}{{\rm dom.dim }}
\newcommand{\E}{{\rm E}}
\newcommand{\Mor}{{\rm Morph}}
\newcommand{\End}{{\rm End}}
\newcommand{\ind}{{\rm ind}}
\newcommand{\rsd}{{\rm res.dim}}
\newcommand{\rd} {{\rm rep.dim}}
\newcommand{\ol}{\overline}
\newcommand{\rad}{{\rm rad}}
\newcommand{\soc}{{\rm soc}}
\renewcommand{\top}{{\rm top}}
\newcommand{\pd}{{\rm proj.dim}}
\newcommand{\id}{{\rm inj.dim}}
\newcommand{\Fac}{{\rm Fac}}
\newcommand{\fd} {{\rm fin.dim }}
\newcommand{\DTr}{{\rm DTr}}
\newcommand{\cpx}[1]{#1^{\bullet}}
\newcommand{\D}[1]{{\mathscr D}(#1)}
\newcommand{\Dz}[1]{{\mathscr D}^+(#1)}
\newcommand{\Df}[1]{{\mathscr D}^-(#1)}
\newcommand{\Db}[1]{{\mathscr D}^b(#1)}
\newcommand{\C}[1]{{\mathscr C}(#1)}
\newcommand{\Cz}[1]{{\mathscr C}^+(#1)}
\newcommand{\Cf}[1]{{\mathscr C}^-(#1)}
\newcommand{\Cb}[1]{{\mathscr C}^b(#1)}
\newcommand{\K}[1]{{\mathscr K}(#1)}
\newcommand{\Kz}[1]{{\mathscr K}^+(#1)}
\newcommand{\Kf}[1]{{\mathscr  K}^-(#1)}
\newcommand{\Kb}[1]{{\mathscr K}^b(#1)}
\newcommand{\modcat}{\ensuremath{\mbox{{\rm -mod}}}}
\newcommand{\Modcat}{\ensuremath{\mbox{{\rm -Mod}}}}
\newcommand{\stmodcat}[1]{#1\mbox{{\rm -{\underline{mod}}}}}
\newcommand{\pmodcat}[1]{#1\mbox{{\rm -proj}}}
\newcommand{\imodcat}[1]{#1\mbox{{\rm -inj}}}
\newcommand{\PI}[1]{#1\mbox{{\rm -prinj}\,}}
\newcommand{\opp}{^{\rm op}}
\newcommand{\otimesL}{\otimes^{\rm\bf L}}
\newcommand{\rHom}{{\rm\bf R}{\rm Hom}}
\newcommand{\projdim}{\pd}
\newcommand{\Hom}{{\rm Hom}}
\newcommand{\Coker}{{\rm Coker}\,\,}
\newcommand{ \Ker  }{{\rm Ker}\,\,}
\newcommand{ \Img  }{{\rm Im}\,\,}
\newcommand{\Ext}{{\rm Ext}}
\newcommand{\StHom}{{\rm \underline{Hom} \, }}
\newcommand{\gm}{{\rm _{\Gamma_M}}}
\newcommand{\gmr}{{\rm _{\Gamma_M^R}}}

\def\vez{\varepsilon}
\def\bz{\bigoplus}
\def\sz {\oplus}
\def\epa{\xrightarrow}
\def\inja{\hookrightarrow}

\newcommand{\lra}{\longrightarrow}
\newcommand{\lraf}[1]{\stackrel{#1}{\lra}}
\newcommand{\ra}{\rightarrow}
\newcommand{\dk}{{\rm dim_{_{\rm k}}}}

{\Large \bf
\begin{center}
Dominant dimensions, derived equivalences and tilting modules
\end{center}}

\begin{center}
\centerline{{\bf Hongxing Chen} and {\bf Changchang Xi$^*$}}
\medskip
School of Mathematical Sciences, BCMIIS, Capital Normal University, \\ 100048 Beijing, People's Republic of  China \\
{\small E-mail: chx19830818@163.com (H.X.Chen), xicc@cnu.edu.cn (C.C.Xi)}\\
\end{center}
\bigskip

\renewcommand{\thefootnote}{\alph{footnote}}
\setcounter{footnote}{-1} \footnote{ $^*$ Corresponding author.
Email: xicc@cnu.edu.cn; Fax: +86 10 68903637.}
\renewcommand{\thefootnote}{\alph{footnote}}
\setcounter{footnote}{-1} \footnote{2010 Mathematics Subject
Classification: 16G10, 18G20; 18E30, 16S50.}
\renewcommand{\thefootnote}{\alph{footnote}}
\setcounter{footnote}{-1} \footnote{Keywords: Derived equivalence; Dominant dimension; Relatively exact sequence; Self-injective algebra; Tilting module.}


\begin{abstract}
The Nakayama conjecture states that an algebra of infinite dominant dimension should  be self-injective. Motivated by understanding this conjecture in the context of derived categories, we study dominant dimensions of algebras under derived equivalences induced by tilting modules, specifically, the infinity of dominant dimensions under tilting procedure. We first give a new method to produce derived equivalences from relatively
exact sequences, and then establish relationships and lower bounds of dominant dimensions for derived equivalences induced by tilting modules. Particularly, we show that under a sufficient condition the infinity of dominant dimensions can be preserved by tilting, and get not only a class of derived equivalences
between two algebras such that one of them is a Morita algebra in the sense of
Kerner-Yamagata and the
other is not, but also the first counterexample to the question whether generalized
symmetric algebras are closed under derived equivalences.
\end{abstract}

\tableofcontents

\section{Introduction}

Derived equivalences play an important role in the representation
theory of algebras and finite groups (\cite{H}, \cite{Broue1994}),
while the Morita theory of derived categories by Rickard provides a
powerful tool to understand these equivalences of rings
\cite{Rickard}. On the one hand, many different approaches to
constructing derived equivalences have been made in recent years. For example, a kind of relatively split sequences has been
introduced in \cite{hx2, hkx} to produce systematically
classical tilting modules and derived equivalences. On the other
hand, many homological invariants of derived equivalences have been
discovered, for instance, Hochschild homology and cohomology,
finiteness of global and finitistic dimensions (see \cite{H, Happel,
keller2, PX}). Unfortunately, so far as we know, there are few papers to
investigate dominant dimensions of algebras in the context of derived
equivalences. Recall that the dominant dimension of an algebra reflects how far the algebra is away from being self-injective, while
the latter forms an important class of algebras in the
representation theory of algebras and has many significant
applications in topology (see \cite{sy} and \cite{koch}). Related to
dominant dimensions, there is a famous open problem: If an algebra
has infinite dominant dimension then it should be self-injective.
This is the so-called Nakayama conjecture (see, for instance, \cite[Conjecture (8), p.410]{ARS}) and has attracted interests of
a lot of mathematicians such as M. Auslander, K. R. Fuller, B. Huisgen-Zimmermann, I. Reiten, H. Tachikawa, and G. V. Wilson.

In the present paper, we shall study behaviors of dominant
dimensions of algebras under derived equivalences, and try to understand the above conjecture in the context of derived categories. For this purpose,
we first have to construct some particular derived equivalences with given dominant dimensions (at least for small dominant dimensions). So,
we introduce a more general notion of relatively exact sequences
which not only capture relatively split sequences defined in
\cite{hx2}, but also provide us with a construction of derived equivalences between
subrings of the endomorphism rings of objects involved. Here, these
subrings can be described explicitly (see Proposition \ref{gder} below),
and the most ideas of the proof of this construction are motivated from \cite{hx2}. However, along the way some of arguments in \cite{hx2}
have been changed and some of the results seem to be new.
As a consequence of these discussions, we get the first main
result, Theorem \ref{main result}, of this paper, which provides a method to construct derived equivalences between algebras with small dominant dimensions.
As a byproduct of this result, we
construct an example of a derived equivalence between a generalized symmetric algebra of dominant dimension $2$
and an algebra of dominant dimension $1$. Thus we answer negatively a question by Ming Fang whether derived equivalences preserve
generalized symmetric algebras.

This example also shows a known phenomenon that, in general, derived
equivalences do not preserve dominant dimensions. So our next purpose is to investigate the relationship of dominant dimensions for derived equivalences induced from tilting modules, and further to consider when such equivalences
preserve dominant dimensions and Morita algebras. In this
direction, our second main result, Theorem \ref{distance}, provides an inequality of the dominant dimension of one algebra in terms of the one of the other, together with the projective dimension of a tilting module. Consequently, we obtain a sufficient condition for tilting procedure to preserve the infinity of dominant dimensions. Moreover, Corollary \ref{equality} gives a lower bound for the dominant dimension of the endomorphism algebra of an arbitrary tilting module $T$ in terms of $T$-gradients of the given algebra, while
Proposition \ref{invariant} provides several sufficient
conditions for tilting procedure to preserve dominant dimensions and Morita algebras. In particular, for a Morita algebra $A$ and any tilting $A$-module $T$, the dominant dimension of $A$ is always less than or equal to the dominant dimension of
the endomorphism algebra $B$ of $T$ plus the projective dimension of $T$; and the endomorphism algebras of canonical tilting $A$-modules
are again Morita algebras and have the same dominant dimension as $A$
does (see Corollaries \ref{compare-Morita} and \ref{cor23}).

The paper is outlined as follows: In Section \ref{sect2}, we fix
notation, introduce the notion of relatively exact sequences with
respect to subcategories, and give a new construction of derived
equivalences. In section \ref{main}, we prove the first main result Theorem \ref{main result}. In section \ref{sect4new}, we prove the second main result Theorem \ref{distance}, where its particular applications to $n$-BB-tilting modules
and canonical tilting modules are also given. In order to describe lower bounds for dominant dimensions of the endomorphism algebras of tilting modules, the notion of $T$-gradients is introduced in this section. The last section, Section \ref{sect4}, is devoted to showing the first
counterexample to a question of whether generalized symmetric
algebras are closed under taking derived equivalences. Also, a few open questions relevant to results in this paper are mentioned there.

\section{Relatively exact sequences and derived equivalences\label{sect2}}

In this section, we first fix some notation, and then generalize
a result in \cite[Thoerem 1.1]{hx2}.

Throughout this section, let $\mathcal C$ be an additive category.

Given two morphisms $f: X\to Y$ and $g: Y\to Z$ in $\mathcal C$, we
denote the composite of $f$ and $g$ by $fg$ which is a morphism from
$X$ to $Z$. The induced morphisms $\Hom_{\mathcal
C}(Z,f):\Hom_{\mathcal C}(Z,X)\ra \Hom_{\mathcal C}(Z,Y)$ and
$\Hom_{\mathcal C}(f,Z): \Hom_{\mathcal C}(Y, Z)\ra \Hom_{\mathcal
C}(X, Z)$ are denoted by $f^*$ and $f_*$, respectively.

Let $X$ be an object in $\mathcal{C}$. Then we denote by $\add(X)$ the full
subcategory of $\mathcal{C}$ consisting of all direct summands of
direct sums of finitely many copies of $X$. The endomorphism algebra
of the object $X$ is denoted by $\End_\mathcal{C}(X)$. It is known
that the Hom-functor $\Hom_\mathcal{C}(X, -)$ is a fully faithful
functor from $\add(X)$ to the category of finitely generated
projective $\End_\mathcal{C}(X))$-modules. We say that $X$ is an
additive generator for $\mathcal C$ if $\add(X)=\mathcal C$.

Let $\C{\mathcal{C}}$ be the category of all complexes over
$\mathcal{C}$ with chain maps, and $\K{\mathcal{C}}$ the homotopy
category of $\C{\mathcal{C}}$. When $\mathcal{C}$ is abelian, the
derived category of $\mathcal{C}$ is denoted by $\D{\mathcal{C}}$,
which is the localization of $\K{\mathcal C}$ at all
quasi-isomorphisms.

Let $R$ be an arbitrary ring with identity. We denote by $R\Modcat$
(respectively, $R\modcat$) the category of all unitary
(respectively, finitely generated) left $R$-modules. As usual, we
simply write $\C{R}$, $\K{R}$ and $\D{R}$ for $\C{R\Modcat}$,
$\K{R\Modcat}$ and $\D{R\Modcat}$, respectively.

Recall that two rings $R$ and $S$ are said to be \emph{derived
equivalent} if $\D{R}$ and $\D{S}$ are equivalent as triangulated
categories. Note that if $R$ and $S$ are finite-dimensional
$k$-algebras over a field $k$, then $R$ and $S$ are derived
equivalent if and only if $\Db{R\modcat}$ and $\Db{S\modcat}$ are equivalent as
triangulated categories. For more details on characterizations of derived equivalences, we refer to \cite{Rickard}, and for some new constructions of derived equivalences, we refer the reader to the recent papers \cite{CYP, hkx, hx2, hx3}.

\medskip
Let us start with recalling from \cite{hx2} the definition of relatively split sequences in
additive categories.

\begin{Def}\label{D-split}
Let $\mathcal{D}$ be a full subcategory of $\mathcal{C}$. A sequence
$$X \lraf{f} M \lraf{g} Y$$
of two morphisms $f$ and $g$ between objects in $\mathcal{C}$ is called a $\mathcal D$-split sequence if

$(1)$ $M\in {\mathcal D};$

$(2)$ $\Hom_\mathcal{C}(D',g)$ and $\Hom_\mathcal{C}(f,D')$ are
surjective for any object $D'\in {\mathcal D};$

$(3)$ $f$ is a kernel of $g$, and $g$ is a cokernel of $f$.
\end{Def}

Given an arbitrary $\cal D$-split sequence, there exists a
derived equivalence between the endomorphism algebras of relevant objects, as shown
by the following result.

\begin{Prop}{\rm\cite[Theorem 1.1]{hx2}}\label{der}
Let $\mathcal{C}$ be an additive category and $M$ an object in
$\mathcal{C}$. Suppose that $$X\lra M' \lra Y $$ is an
$\add(M)$-split sequence in $\mathcal{C}$.  Then
$\End_\mathcal{C}(X\oplus M)$ and $\End_\mathcal{C}(Y\oplus M)$ are
derived equivalent.
\end{Prop}

For our purpose, we introduce the following definition of $\cal D$-exact sequences, which modifies and generalizes slightly the one of $\mathcal D$-spit sequences.

\begin{Def}\label{ass}
Let $\mathcal{D}$ be a full subcategory of $\mathcal C$ . A sequence
$$X \lraf{f} M_0 \lraf{g} Y$$
of objects and morphisms in $\mathcal{C}$ is called a $\mathcal D$-exact
sequence provided that

$(1)$ $M_0\in{\mathcal D}$.

$(2)$ The following two sequences of abelian groups are exact:
$$(\dag)\quad\;
0\lra \Hom_\mathcal{C}(X\oplus M, \,X)\lraf {f^*}
\Hom_\mathcal{C}(X\oplus M,\, M_0)\lraf{g^*}
\Hom_\mathcal{C}(X\oplus M, \,Y)$$
$$(\ddag)\quad\;
0\lra \Hom_\mathcal{C}(Y, \,M\oplus Y)\lraf {g_*}
\Hom_\mathcal{C}(M_0, \,M\oplus Y)\lraf{f_*} \Hom_\mathcal{C}(X,
\,M\oplus Y)$$ for every object $M$ in $\mathcal D$.
\end{Def}

Note that the condition $(2)$ in Definition \ref{ass} implies
$fg=0$. Moreover, if $f$ is a kernel of $g$ and $g$ is a cokernel of
$f$, then the condition (2) holds automatically. Thus $\mathcal
D$-split sequences in $\mathcal C$ are $\mathcal D$-exact
sequences in $\mathcal C$. But the converse is not true: Since every short exact sequence $0\ra X\ra M\ra Y\ra 0$ in an abelian  category is an $\add(M)$-exact sequence, we get not only the ubiquity of relatively exact sequences, but also examples of $\mathcal D$-exact sequences which is not $\mathcal D$-split. For instance, we take $A=k[T_1,T_2]/(T_1^2,T_2^2, T_1T_2)$ with $k$ a field, $\mathcal{C}=A\modcat$ and $X=A/\rad(A)$, then there is a short exact sequence
$ 0\ra X\ra A\ra A/(T_1)\ra 0$, for which the condition (2) in Definition \ref{D-split} is not satisfied. Thus this sequence is not $\add(_AA)$-split.

\medskip
In the following, we shall focus on the most interesting case where $\mathcal{ D}=
\add(M)$ for $M$ an object in $\mathcal{C}$. Observe that, for an arbitrary
$\add(M)$-exact sequence (not necessarily an $\add(M)$-split
sequence), we do not have to get a derived equivalence between endomorphism
rings as in Proposition \ref{der}. However, we shall prove that there does exist a derived
equivalence between subrings of corresponding endomorphism rings.

\begin{Prop}\label{gder}
Let $\mathcal{C}$ be an additive category and $M$ an object in
$\mathcal{C}$. Suppose $$X\lraf{f} M_0 \lraf{g} Y $$ is an
$\add(M)$-exact sequence in $\mathcal{C}$. Set
$$R:=\left\{
{\small \left(\begin{array}{lc} \;\,h_1& h_2\\
 f h_3 & h_4 \end{array}\right)\in
\left(\begin{array}{lc} \;\End_\mathcal{C}(M) &\Hom_\mathcal{C}(M, X)\\
\Hom_\mathcal{C}(X, M)& \End_\mathcal{C}(X) \end{array}\right)}
\bigg |
\begin{array}{ll}
h_3\in\Hom_\mathcal{C}(M_0, M) \;\mbox{\rm and there exists} \\
h_5 \in\End_\mathcal{C}(M_0)\; \mbox{\rm such that\;} h_4 f= f h_5
\end{array}
\right\}
$$
and
$$S:=\left\{
{\small \left(\begin{array}{lc} \;h_1& h_2\, g\\
 h_3 & h_4 \end{array}\right)\in
\left(\begin{array}{lc} \;\End_\mathcal{C}(M) &\Hom_\mathcal{C}(M, Y)\\
\Hom_\mathcal{C}(Y, M)& \End_\mathcal{C}(Y) \end{array}\right)}
\bigg |
\begin{array}{ll}
h_2\in\Hom_\mathcal{C}(M, M_0) \;\mbox{\rm and there exists} \\
h_5 \in\End_\mathcal{C}(M_0)\; \mbox{\rm such that\;} g h_4=h_5 g
\end{array}
\right\}.
$$

\smallskip
\noindent Then $R$ and $S$ are subrings of $\End_\mathcal{C}(M\oplus
X)$ and $\End_\mathcal{C}(M\oplus Y)$, respectively. Moreover, they
are derived equivalent.
\end{Prop}

{\it Proof.} The proof here is actually motivated by \cite[Lemma
3.4]{hx2}.

Set $V:=X\oplus M$ and $\Lambda:=\End_\mathcal{C}(V)$. By the exact
sequence $(\dag)$ in Definition \ref{ass}, there exists an exact
sequence
$$0\lra \Hom_\mathcal{C}(V, \,X)\lraf {f^*}
\Hom_\mathcal{C}(V,\, M_0)\lraf{g^*} \Hom_\mathcal{C}(V, \,Y)$$ of
$\Lambda$-modules. Further, we define $L:=\Img(g^*)$, the image of the
map $g^*$, and $T:=\Hom_\mathcal{C}(V, M)\oplus L$.

In the following, we shall divide the whole proof of Proposition
\ref{gder} into four steps.

\smallskip
$(1)$ We claim that $\End_{\Lambda}(T)\simeq S$ as rings.

To show this, let $\ol{f}=(0, f): X \ra M\oplus M_0$ and $\ol{g}=\left(\begin{array}{ll} 1 &0\\
0 & g\end{array}\right): M\oplus M_0\ra M\oplus Y.$ Then, from $f\,g=0$
we have $\ol{f}\,\ol{g}=0$. Moreover, there exists the following
exact sequence of $\Lambda$-modules:
$$0\lra \Hom_\mathcal{C}(V, \,X)\lraf {{\ol f\,}^*}
\Hom_\mathcal{C}(V,\,M\oplus M_0)\lraf{{\ol g\,}^*} T\lra 0.
$$
Let $\cpx{P}$ be the following complex: $$ 0\lra \Hom_\mathcal{C}(V,
\,X)\lraf {{\ol f\,}^*} \Hom_\mathcal{C}(V,\,M\oplus M_0)\lra 0$$
with $\Hom_\mathcal{C}(V, \,X)$ in degree $-1$. Note
that both $\Hom_\mathcal{C}(V, X)$ and $\Hom_\mathcal{C}(V,\,M\oplus
M_0)$ are finitely generated projective $\Lambda$-modules
since $X\in\add(V)$ and $M_0\in\add(M)\subseteq \add(V)$. This
implies that $\End_{\Lambda}(T)\simeq \End_{\D
{\Lambda}}(\cpx{P})\simeq \End_{\K {\Lambda}}(\cpx{P})$ as rings. As
the Hom-functor $\Hom_\mathcal{C}(V, -):\add(V)\to
\add(_\Lambda\Lambda)$ is fully faithful, we see that $\End_{\K
{\Lambda}}(\cpx{P})\simeq \End_{\K{\add(V)}}(\cpx{Q})$ as rings,
where $\cpx{Q}$ is defined to be the complex: $$ 0\lra X \lraf{\ol
f} M\oplus M_0\lra 0$$ where the object $X$ is of degree $-1$.
To finish the proof of (1), it suffices to show that there exists an injective ring homomorphism
$\Psi: \End_{\K{\add(V)}}(\cpx{Q})\ra \End_\mathcal{C}(M\oplus Y)$
such that $\Img(\Psi)=S$ because it then follows that
$\End_{\Lambda}(T)\simeq\End_{\K{\add(V)}}(\cpx{Q})\simeq S$ as
rings, as claimed.

Let $(\alpha, \beta):\cpx{Q}\to \cpx{Q}$ be an arbitrary chain map
with $\alpha\in \End_\mathcal{C}(X)$ and
$\beta\in\End_\mathcal{C}(M\oplus M_0)$. Then $\alpha
\ol{f}=\ol{f}\beta$. Now, we point out that there exists a unique
morphism $\gamma\in\End_\mathcal{C}(M\oplus Y)$ such that the
following diagram is commutative in $\mathcal{C}$:
$$
\xymatrix {X\ar[d]_-{\alpha} \ar[r]^-{\ol f}& M\oplus
M_0\ar[d]_-{\beta}\ar[r]^-{\ol g}& M\oplus Y\ar@{-->}[d]^-{\gamma}\\
X\ar[r]^-{\ol f}& M\oplus M_0 \ar[r]^-{\ol g}& M\oplus Y.}
$$
Actually, by the sequence $(\ddag)$ in Definition \ref{ass} , we
have the following exact sequence of abelian groups:
$$0\lra \Hom_\mathcal{C}(M\oplus Y, \,M\oplus Y)\lraf {{\ol g\,}_*}
\Hom_\mathcal{C}(M\oplus M_0, \,M\oplus Y)\lraf{{\ol f\,}_*}
\Hom_\mathcal{C}(X, \,M\oplus Y).$$ Since ${\ol
f\,}_*(\beta\,\ol{g})=\ol{f}\,\beta\,\ol{g}=\alpha\,\ol{f}\,\ol{g}=0$,
there is a unique morphism $\gamma\in\End_\mathcal{C}(M\oplus Y)$
such that ${\ol g\,}_*(\gamma)=\ol{g}\,\gamma=\beta\,\ol{g}$.

Now, we prove that the chain map $(\alpha, \beta)$ is homotopic to the
zero map if and only if $\gamma=0$.

In fact, if $(\alpha, \beta)$ is null-homotopic, then there exists a morphism
$\delta: M\oplus M_0 \to X$ such that $\alpha=\ol{f}\,\delta$ and
$\beta=\delta\,\ol{f}$. In this case, we have
$\beta\,\ol{g}=\delta\,\ol{f}\,\ol{g}=0$, and therefore $\gamma=0$.

Suppose that $\gamma=0$. Then $0=\ol{g}\,\gamma=\beta\,\ol{g}$.
Since $M\oplus M_0\in\add(M)$, we know from the sequence $(\dag)$ in
Definition \ref{ass} that the following sequence  $$0\lra
\Hom_\mathcal{C}(M\oplus M_0, \,X)\lraf {{\ol f\,}^*}
\Hom_\mathcal{C}(M\oplus M_0,\, M\oplus M_0)\lraf{{\ol g\,}^*}
\Hom_\mathcal{C}(M\oplus M_0, M\oplus Y)$$ is exact. This implies
that there exists a morphism $\sigma: M\oplus M_0 \to X$ such that
$\beta=\sigma\,\ol{f}$. Consequently, we have
$(\alpha-\ol{f}\,\sigma)\,\ol{f}=0$ due to $\alpha
\ol{f}=\ol{f}\beta$. Similarly, by the sequence $(\dag)$ in
Definition \ref{ass}, we see that the induced map
$\Hom_\mathcal{C}(X,\ol{f}): \Hom_\mathcal{C}(X, X)\ra
\Hom_\mathcal{C}(X, M\oplus M_0)$ is injective. This gives rise to
$\alpha-\ol{f}\,\sigma=0$, that is $\alpha=\ol{f}\,\sigma$. Thus
$(\alpha, \beta)$ is null-homotopic.

As a result, the following map
$$\Psi:\; \End_{\K{\add(V)}}(\cpx{Q})\lra
\End_\mathcal{C}(M\oplus Y),\;\, \ol{(\alpha, \beta)}\mapsto
\gamma$$ is well defined. Clearly, this map is an injective ring
homomorphism. It remains to show $\Img(\Psi)=S$.

To check this equality, we write $$\beta={\small \left(\begin{array}{ll} \beta_1 &\beta_2\\
\beta_3 & \beta_4\end{array}\right)\in
\left(\begin{array}{lc} \;\End_\mathcal{C}(M) &\Hom_\mathcal{C}(M, M_0)\\
\Hom_\mathcal{C}(M_0, M)& \End_\mathcal{C}(M_0)
\end{array}\right)}, \;
\gamma={\small \left(\begin{array}{ll} \gamma_1 &\gamma_2\\
\gamma_3 & \gamma_4\end{array}\right)\in
\left(\begin{array}{lc} \;\End_\mathcal{C}(M) &\Hom_\mathcal{C}(M, Y)\\
\Hom_\mathcal{C}(Y, M)& \End_\mathcal{C}(Y) \end{array}\right).}$$
From $\ol{g}\,\gamma=\beta\,\ol{g}$, we obtain
$\gamma_1=\beta_1$, $\gamma_2=\beta_2\,g$, $\,\beta_3=g\,\gamma_3$
and $g\,\gamma_4=\beta_4\,g$. This implies that $\gamma\in S$ and
$\Img(\Psi)\subseteq S$.

Conversely, suppose $$h:=
{\small \left(\begin{array}{lc} \;h_1& h_2\, g\\
 h_3 & h_4 \end{array}\right)\in S}$$where
$h_2\in \Hom_\mathcal{C}(M, M_0)$ and  there is a morphism $h_5
\in\End_\mathcal{C}(M_0)$ such that $g h_4=h_5 g.$ Define
$\beta:={\small \left(\begin{array}{ll}\; h_1 &h_2\\
g\,h_3 & h_5\end{array}\right)}\in\End_\mathcal{C}(M\oplus M_0).$
Then one can verify that $\ol{g}\,h=\beta\,\ol{g}.$ Since
$\ol{f}\,\beta\,\ol{g}=\ol{f}\,\ol{g}\,h=0$, we know from the
sequence $(\dag)$ in Definition \ref{ass} that there exists a unique
morphism $\alpha\in\End_\mathcal{C}(X)$ such that the following
diagram commutes in $\mathcal{C}$:
$$
\xymatrix {X\ar@{-->}[d]_-{\alpha} \ar[r]^-{\ol f}& M\oplus
M_0\ar[d]_-{\beta}\ar[r]^-{\ol g}& M\oplus Y\ar[d]^-{h}\\
X\ar[r]^-{\ol f}& M\oplus M_0 \ar[r]^-{\ol g}& M\oplus Y.}
$$
This implies that $\Psi\big(\ol{(\alpha, \beta)}\big)=h$ and
$S\subseteq \Img(\Psi)$. Thus $S=\Img(\Psi)$. Since $\Psi$ is an
injective ring homomorphism, we infer that $S$ is actually a subring
of $\End_\mathcal{C}(M\oplus Y)$ and that $\Psi:
\End_{\K{\add(V)}}(\cpx{Q})\to S$ is an isomorphism of rings.

Hence $\End_{\Lambda}(T)\simeq \End_{\K{\add(V)}}(\cpx{Q})\simeq S$
as rings. This finishes the proof of $(1)$.

\smallskip
$(2)$ We claim that if the induced map $\Hom_\mathcal{C}(f, M):
\Hom_\mathcal{C}(M_0, M)\ra \Hom_\mathcal{C}(X, M)$ is surjective,
then $\Lambda$ is derived equivalent to $S$.

Set $N:=\Hom_\mathcal{C}(V, M)$. Then $_\Lambda T=N\oplus L$ by the
foregoing notation, where the module $L$ arises in the following
exact sequence of $\Lambda$-modules:
$$(\ast)\quad
0\lra \Hom_\mathcal{C}(V, \,X)\lraf { f^*} \Hom_\mathcal{C}(V,\,
M_0)\lraf{g^*} L\lra 0.
$$
Clearly, we have $\Hom_\mathcal{C}(V,\, M_0)\in\add(N)$ due to
$M_0\in\add(M)$. On the one hand, since $_{\Lambda}N=\Hom_\mathcal{C}(V, M)$ is
a finitely generated projective $\Lambda$-module, the induced map $\Hom_\Lambda(N,
g^*):\Hom_\Lambda(N, \Hom_\mathcal{C}(V, M_0))\lra \Hom_\Lambda(N,
L)$ is naturally surjective. On the other hand, since the functor
$\Hom_\mathcal{C}(V, -):\add(V)\to \add(_\Lambda\Lambda)$ is fully
faithful, we know that the induced map $$\Hom_\Lambda(f^*, N):
\Hom_\Lambda(\Hom_\mathcal{C}(V, M_0), N)\lra
\Hom_\Lambda(\Hom_\mathcal{C}(V, X), N)$$ is surjective if and only
if so is the map $\,\Hom_\mathcal{C}(f, M): \Hom_\mathcal{C}(M_0,M)\ra \Hom_\mathcal{C}(X,M)$.

Assume that $\Hom_\mathcal{C}(f, M)$ is surjective. By Definition
\ref{D-split}, the exact sequence $(\ast)$ is an $\add(_\Lambda
N)$-split sequence. Now, it follows from Proposition \ref{der} that the
rings $\Lambda$ and $\End_\Lambda(T)$ are derived equivalent. Since
$\End_\Lambda(T)$ is isomorphic to $S$ by $(1)$, we know that
$\Lambda$ is derived equivalent to $S$. This finishes the proof of
$(2)$.

\smallskip
$(3)$ Since the notion of $\add(M)$-exact sequences (see
Definition \ref{ass}) is self-dual, we can show
dually that if the induced map $\Hom_\mathcal{C}(M, g):
\Hom_\mathcal{C}(M, M_0)\ra \Hom_\mathcal{C}(M, Y)$ is surjective,
then $\End_\mathcal{C}(M\oplus Y)$ is derived equivalent to $R$.

\smallskip
$(4)$ Now we show that $R$ and $S$ are derived equivalent.

Recall that the sequence $(\ast)$ is an exact sequence of
$\Lambda$-modules. Certainly, it is an $\add(_\Lambda N)$-exact
sequence in the category of $\Lambda$-modules. Moreover, this
sequence always has the following property: The induced map
$\Hom_\Lambda(N, g^*):\Hom_\Lambda(N, \Hom_\mathcal{C}(V, M_0))\ra
\Hom_\Lambda(N, L)$ is surjective. Set $U:=\Hom_\mathcal{C}(V, X)$.
Applying $(3)$ to the sequence ($*$), we see that $\End_\Lambda(T)$ is derived
equivalent to the ring
$$\widetilde{R}:=\left\{
{\small \left(\!\begin{array}{lc} \;\,h_1& h_2\\
 f h_3 & h_4 \end{array}\!\right)\in
\left(\!\begin{array}{lc} \;\End_\Lambda(N) &\Hom_\Lambda(N, U)\\
\Hom_\Lambda(U, N)& \End_\Lambda(U)
\end{array}\!\right)} \bigg |
\begin{array}{ll}
h_3\in\Hom_\Lambda(\Hom_\mathcal{C}(V, M_0), N), \;\mbox{\rm there
exists} \\ h_5 \in\End_\Lambda(\Hom_\mathcal{C}(V, M_0))\;
\mbox{\rm such that\;} h_4 f^*= f^* h_5
\end{array}
\right\}.
$$
Since the functor $\Hom_\mathcal{C}(V, -):\add(V)\to
\add(_\Lambda\Lambda)$ is fully faithful, one can easily check that
this functor induces a ring isomorphism from $R$ to $\widetilde{R}$.
Recall from $(1)$ that $\End_\Lambda(T)$ is
isomorphic to $S$. Thus $R$ and $S$ are derived equivalent. This
finishes the proof of Proposition \ref{gder}. $\square$

\medskip
Remarks. (1) If $X\ra M_0\ra Y$ is an $\add(M)$-split sequence in $\mathcal C$, then $R=\End_{\mathcal C}(M\oplus X)$ and $S=\End_{\mathcal C}(M\oplus Y)$. Thus Proposition \ref{gder} implies Proposition \ref{D-split}.

(2) If $\mathcal{C}$ is an abelian category and the sequence in Proposition \ref{der} is a short exact sequence in $\mathcal{C}$, then the derived equivalence between $R$ and $S$ follows also from \cite[Corollary 3.4]{CYP}.

As an easy consequence of the above proposition, we know that the rings $R$ and $S$ have the same algebraic $K$-groups since derived equivalences preserve algebraic $K$-theory, and that the finitistic dimension of $R$ is finite if and only if so is the one of $S$ (see \cite{PX}). Further applications will be given in the next section.

\section{Dominant dimensions and derived equivalences\label{main}}
Throughout this section, $k$ stands for a fixed field. All algebras
considered are finite-dimensional $k$-algebras with identity, and all modules are
finitely generated left modules.

\subsection{Basic facts on dominant dimensions}

Let $A$ be an algebra. We denote by $\rad(A)$ the Jacobson radical
of $A$, by $A\pmodcat$ (respectively, $A\imodcat$) the full
subcategory of $A\modcat$ consisting of projective (respectively,
injective) modules,  by $D$ the usual $k$-duality $\Hom_{k}(-, k)$,
and by $\nu_{A}$ the Nakayama functor $D\Hom_{A}(-,\,_{A}A)$ of $A$.
Note that $\nu_A$ is an equivalence from $A\pmodcat$ to $A\imodcat$
with the inverse $\nu^-_A=\Hom_A(D(A),-)$.
The category of projective-injective $A$-modules is denoted by $A\PI$.

Let $X$ be an $A$-module. By $\Omega^i_A(X)$, $\soc(X)$ and $I(X)$
we denote the $i$-th syzygy for $i\in\mathbb{Z}$, the socle and the
injective envelope of $X$, respectively.

\medskip
For an $A$-module $X$, we consider its minimal injective resolution
$$0\lra {}_AX\lra I_0\lra I_1\lra I_2\lra\cdots. $$
Let $I$ be an injective $A$-module and $0\le n\le \infty$. If $n$ is maximal with the property
that all modules $I_j$ are in $\add(I)$ for $j < n$, then $n$ is called the \emph{dominant dimension
of $X$ with respect to $I$}, denoted by $I$-$\dm(X)$. If $\add(I)$ = $A\PI$, we simply
write $\dm(X)$ and call it the \emph{dominant dimension} of $X$. Since $\dm(_AA) =\dm(A_A)$ by \cite[Theorem 4]{Muller}, we just write $\dm(A)$ and call it the \emph{dominant dimension
of $A$}. It is clear that $\dm(A) = \min\{\dm(P)\mid P\in \add(_AA)\}$.

If $A$ is self-injective, that is, the regular module $_AA$ is injective, then $\dm(A)=\infty$. The converse of this statement is the well-known, longstanding Nakayama conjecture: If $\dm(A)=\infty$, then $A$ is self-injective. Equivalently, if $A$ is not self-injective, then $\dm(A)\le m$ for a positive integer $m$. Hence, in order to understand this conjecture, it makes sense to investigate upper bounds for dominant dimensions.

It is well known that $\dm(A)\geq 2$ if and only if there
exists an algebra $B$ and a generator-cogenerator $V$ over $B$ such
that $A\simeq \End_B(V)$ as algebras (see \cite[Theorem 2]{Muller}).
In fact, let $e$ be an idempotent element of $A$ such that
$\add(\nu_A(Ae))$ coincides with the full subcategory of $A\modcat$
consisting of all projective-injective $A$-modules. If $\dm(A)\geq
2$, then we can choose $B=eAe$ and $V=eA$. Furthermore, in this
case, $\dm(A)=n$ if and only if $\Ext_{eAe}^i(eA, eA)=0$ for all
$1\leq i\leq n-2$ and $\Ext_{eAe}^{n-1}(eA, eA)\neq 0$ (see
\cite[Lemma 3]{Muller}). Thus, for a self-injective algebra $A$ and  an $A$-module $Y$ without projective summands,
if there is an $n\ge 0$ such that $\Ext^{n+1}_A(Y,Y)\ne 0$ and $\Ext^j_A(Y,Y)=0$ for all $1\le j\le n$, then $\dm (\End_A(A\oplus Y))= n+2.$

\medskip
To estimate dominant dimensions, we need the following result which is essentially taken from \cite[Lemma 1.1 and Corollary 1.3]{Miyachi}.

\begin{Lem}\label{inj-resol}
Let $\Lambda$ be an algebra and let
$0\to Y_{-1}\to Y_0\to Y_1\to Y_2\to \cdots\to Y_{m-1}\to Y_m\to 0$ be a long exact sequence of $\Lambda$-modules with $m\geq 0$.  Then the following statements are true:

$(1)$ If each $Y_i$ with $0\leq i\leq m$ has an injective resolution $0\to Y_i\to I_i^0\to I_i^1\to I_i^2\to\cdots$, then $Y_{-1}$ has an injective resolution of the following form:
$$
0\lra Y_{-1}\lra I_0^0\lra\bigoplus_{0\leq r\leq \min\{m,\,1\}}I_r^{1-r} \lra \bigoplus_{0\leq r\leq \min\{m,2\}}I_r^{2-r}\lra \cdots\lra \bigoplus_{0\leq r\leq \min\{m,\,s\}}I_r^{s-r}\lra\cdots
$$

$(2)$ If each $Y_j$ with $-1\leq j\leq m-1$ has an injective resolution $0\to Y_j\to I_j^0\to I_j^1\to I_j^2\to\cdots$, then $Y_m$ has an injective resolution of the following form:
$$
0\lra Y_{m}\lra Q_{m-1}\lra \bigoplus_{-1\leq r\leq m-1}I_r^{\,m-r} \lra\bigoplus_{-1\leq r\leq m-1}I_r^{\,m+1-r}\lra\cdots\lra\bigoplus_{-1\leq r\leq m-1}I_r^{\,m+s-r}\lra\cdots
$$
where $Q_{m-1}$ is a direct summand of the module $\bigoplus_{r=-1}^{m-1}I_r^{\,m-1-r} $.
\end{Lem}

As a consequence of Lemma \ref{inj-resol}, we have the following result.

\begin{Koro}\label{calculating dd}
Let $\Lambda$ be an algebra and let $I$ be an injective $\Lambda$-module. Suppose that
$0\to Y_{-1}\to Y_0\to Y_1\to Y_2\to \cdots\to Y_{m-1}\to Y_m\to 0$ is a long exact sequence of $\Lambda$-modules with $m\geq 0$. Then

$(1)$ $I\emph{-}\dm(Y_{-1})\geq \min\big\{I\emph{-}\dm(Y_j)+j\mid 0\leq j\leq m\big\}$.

$(2)$ $I\emph{-}\dm(Y_m)\geq \min\big\{I\emph{-}\dm(Y_j)+j\mid -1\leq j\leq m-1\big\}-m+1$.
\end{Koro}

{\it Proof.} $(1)$ Let $t:=\min\big\{I\emph{-}\dm(Y_j)+j\mid 0\leq j\leq m\big\}$. Then  $I\emph{-}\dm(Y_j)\geq t-j$ for each $0\leq j\leq m$. For such $j$, let
$$
0\lra Y_j\lra I_j^0\lra I_j^1\lra I_j^2\lra \cdots
$$
be a minimal injective resolution of $Y_j$. Then $I_j^u\in\add(I)$ for each
$0\leq u\leq t-j-1$. This implies that
$I_r^{s-r}\in\add(I)$ for $0\leq r\leq \min\{m,s\}$ and $s\leq t-1$.
By Lemma \ref{inj-resol}(1), the module $Y_{-1}$ has an injective resolution
$$
0\lra Y_{-1}\lra I^0\lra I^1\lra I^2\lra \cdots\lra I^{t-1}\lra I^{t}\lra \cdots
$$
such that $I^{j}\in\add(I)$ for all $0\leq j\leq t-1$. Thus
$I\emph{-}\dm(Y_{-1})\geq t$.

$(2)$ Let $m'=\min\big\{I\emph{-}\dm(Y_j)+j\mid -1\leq j\leq m-1\big\}$.
If $m'\leq m-1$, then $I\emph{-}\dm(Y_m)\geq 0\geq m'-m+1$ and therefore $(2)$ holds. Now, we suppose that $m'\geq m$. Then $I\emph{-}\dm(Y_j)\geq m'-j\geq 1$ for $-1\leq j\leq m-1$.
For each $j$, let
$$ 0\lra Y_j\lra I_j^0\lra I_j^1\lra I_j^2\lra \cdots$$
be a minimal injective coresolution of $Y_j$.  Then $I_j^{p}\in\add(I)$ for $0\leq p\leq m'-j-1$.
This implies that $I_j^{q-j}\in\add(I)$ for all $-1\leq j\leq m-1\leq q\leq m'-1$. By Lemma \ref{inj-resol}(2), the module $Y_{m}$ has an injective resolution
$$
0\lra Y_{m}\lra E^0\lra E^1\lra E^2\lra \cdots\lra E^{i-1}\lra E^{i}\lra \cdots
$$
such that $E^{i}\in\add(I)$ for all $0\leq i\leq m'-m$. Consequently, we obtain
$I\emph{-}\dm(Y_{m})\geq m'-m+1$. $\square$

\medskip
Recall that a homomorphism $f: Y_0\to X$ of $A$-modules is called a minimal right $\add(Y)$-approximation of $X$ if $f$ is minimal, $Y_0\in\add(Y)$ and the map $\Hom_A(Y,Y_0)\to \Hom_A(Y,X)$ is surjective. A complex of $A$-modules of the following form
$$
 \cdots\lra Y^{-i}\lraf{f^{-i}} Y^{-i+1}\lraf{f^{-i+1}} \cdots\lraf{f^{-3}} Y^{-2}\lraf{f^{-2}} Y^{-1}\lraf{f^{-1}} X\lraf{f^0} 0$$
with $X$ in degree $0$, is called a \emph{minimal right $\add(Y)$-approximation sequence} of $X$ if the homomorphism $Y^{-i}\to \Ker(f^{-i+1})$ induced from $f^{-i}$ is a minimal right $\add(Y)$-approximation of $\Ker(f^{-i+1})$ for each $i\geq 1$.
Note that, up to isomorphism of complexes, such a sequence is unique and depends only on $X$.

\medskip
The following result is useful for calculation of relative dominant dimensions of modules.

\begin{Lem}\label{idempotent}
Let $P=Ae$ and $I=\nu_A(Ae)$ with $e^2=e\in A$. For any $A$-module $X$, we have
$$I\mbox{-}\dm(X)=\inf\{i\geq 0\mid \Hom_{\K{A}}(\cpx{P}, X[i])\neq 0\},$$
where the complex
$$\cpx{P}:\quad \cdots\lra P^{-i}\lra P^{-i+1}\lra \cdots\lra P^{-2}\lra P^{-1}\lra {_A}A\lra 0$$ is a minimal right $\add(_AP)$-approximation sequence of $_AA$.
\end{Lem}

{\it Proof.}
For each $n\in\mathbb{N}$, it follows from \cite[Proposition 2.6]{APT} that $I$-$\dm(X)\geq n+1$ if and only if $\Ext_A^i(A/AeA,X)=0$ for all $0\leq i\leq n$. The latter is also equivalent to $\Ext_A^i(Y,X)=0$ for all $A/AeA$-modules $Y$ and for all $0\leq i\leq n$. Thus
$$I\emph{-}\dm(X)=\inf\{i\geq 0\mid\Ext_A^i(A/AeA,X)\neq 0\}.$$
To show Lemma \ref{idempotent}, it suffices to check that
$$
\inf\{i\geq 0\mid\Ext_A^i(A/AeA,X)\neq 0\}=\inf\{i\geq 0\mid \Hom_{\K{A}}(\cpx{P}, X[i])\neq 0\}.
$$
As a preparation, we first prove the following statement:

\medskip
$(\ast)\;$
Let $\cpx{Q}:=(Q^j)_{j\in\mathbb{Z}}$ be a complex of $A$-modules such that $Q^j=0$ for $j>0$
and that $H^j(\cpx{Q})\in A/AeA\modcat$ for $j\leq 0$. If $\Ext_A^i(A/AeA,X)=0$ for all
$0\leq i\leq n-1$, then $\Hom_{\D{A}}(\cpx{Q},X[m])=0$ for $m\leq n-1$ and $\Hom_{\D{A}}(\cpx{Q},X[n])\simeq \Hom_{\D{A}}(H^0(\cpx{Q}),X[n])=0$.

\medskip
To prove this statement, we suppose that $\Ext_A^i(A/AeA,X)=0$ for all $0\leq i\leq n-1$. Then $\Hom_{\D{A}}(H^j(\cpx{Q}),X[i])\simeq\Ext_A^i(H^j(\cpx{Q}),X)=0$ for all $j\leq 0$. Note that $\Hom_{\D{A}}(H^j(\cpx{Q}),X[r])=0$ for any $r\leq 0$. Thus $\Hom_{\D{A}}(H^j(\cpx{Q}),X[s])=0$ for all $s\leq n-1$. Now, we take iterated canonical truncations of complexes, and obtain a series of distinguish triangles in $\D{A}$:
$$
\tau^{\leq j-1}(\cpx{Q})\lra\tau^{\leq j}(\cpx{Q})\lra H^j(\cpx{Q})[-j]\lra \tau^{\leq j-1}(\cpx{Q})[1]
$$
where $j\leq 0$ and $\tau^{\leq 0}(\cpx{Q})=\cpx{Q}$. It follows that, for any $m\leq n-1$, we have $$\Hom_{\D{A}}(\tau^{\leq j}(\cpx{Q}),X[m])\hookrightarrow\Hom_{\D{A}}(\tau^{\leq j-1}(\cpx{Q}),X[m]).$$ This leads to the following inclusions of abelian groups:
{\footnotesize $$\Hom_{\D{A}}(\cpx{Q},X[m])\hookrightarrow\Hom_{\D{A}}(\tau^{\leq -1}(\cpx{Q}),X[m])\hookrightarrow\cdots\hookrightarrow\Hom_{\D{A}}(\tau^{\leq -m}(\cpx{Q}),X[m])\hookrightarrow\Hom_{\D{A}}(\tau^{\leq -m-1}(\cpx{Q}),X[m]).$$}
Since $\tau^{\leq -m-1}(\cpx{Q})$ can only have non-zero terms in degrees smaller than $-m$, we get $\Hom_{\D{A}}(\tau^{\leq -m-1}(\cpx{Q}),X[m])=0$, and therefore $\Hom_{\D{A}}(\cpx{Q},X[m])=0$. Since $\big(\tau^{\leq -1}(\cpx{Q})[-1]\big)^t=0$ for any $t\geq 0$ and $$H^t\big(\tau^{\leq -1}(\cpx{Q})[-1]\big)\simeq H^{t-1}(\cpx{Q})\in A/AeA\modcat$$ for any $t\leq 0$, we know that $\Hom_{\D{A}}(\tau^{\leq -1}(\cpx{Q}),X[n])\simeq \Hom_{\D{A}}(\tau^{\leq -1}(\cpx{Q})[-1],X[n-1])=0$ and that  $\Hom_{\D{A}}(\tau^{\leq -1}(\cpx{Q})[1],X[n])\simeq \Hom_{\D{A}}(\tau^{\leq -1}(\cpx{Q})[-1],X[n-2])=0$.  So if we apply $\Hom_{\D{A}}(-,X[n])$ to the triangle $\tau^{\leq -1}(\cpx{Q})\to \cpx{Q}\to H^0(\cpx{Q})\to \tau^{\leq -1}(\cpx{Q})[1]$,
then $\Hom_{\D{A}}(\cpx{Q},X[n])\simeq \Hom_{\D{A}}(H^0(\cpx{Q}),X[n]).$
This finishes the proof of $(\ast)$.

Next, we show that $\Ext_A^i(A/AeA,X)=0$ for all $0\leq i\leq n$ if and only if $\Hom_{\K{A}}(\cpx{P},X[m])=0$ for all $0\leq m\leq n$. This leads to $\inf\{i\geq 0\mid\Ext_A^i(A/AeA,X)\neq 0\}=\inf\{i\geq 0\mid \Hom_{\K A}(\cpx{P}, X[i])\neq 0\}.$

In fact, since $\cpx{P}$ is an above-bounded complex of projective $A$-modules, we have
$\Hom_{\K{A}}(\cpx{P},X[m])\simeq \Hom_{\D{A}}(\cpx{P},X[m]).$
So, it suffices to show that $\Ext_A^i(A/AeA,X)=0$ for all $0\leq i\leq n$ if and only if $\Hom_{\D{A}}(\cpx{P},X[m])=0$ for all $0\leq m\leq n$.

By assumption, the complex $\cpx{P}$ is a minimal right $\add(P)$-approximation sequence of $_AA$. It follows that the complex $\cpx{\Hom}_A(P,\cpx{P})$ is exact. Since $P:=Ae$ is projective, we see that $eH^j(\cpx{P})=0$ for any $j\leq 0$. In other words, $H^j(\cpx{P})\in A/AeA\modcat$. As the $(-1)$-th differential $f^{-1}: P^{-1}\to A$ in $\cpx{P}$ is a minimal right $\add(P)$-approximation of $_AA$, we have $\Img(f^{-1})=AeA$, and therefore $H^0(\cpx{P})=\Coker(f^{-1})=A/AeA$. This gives rise to $\Hom_A(A/AeA,X)=\Hom_A(H^0(\cpx{P}),X)\simeq \Hom_{\D{A}}(\cpx{P},X)$. Thus $\Hom_A(A/AeA,X)=0$ if and only if $\Hom_{\D{A}}(\cpx{P},X)=0$. Now, with the help of the fact $(\ast)$, one can verify the above statement by induction on $n$. $\square$

\medskip
A projective $A$-module $P$ is said to be \emph{$\nu$-stable} if $\nu_A^i(P)$ are projective for all $i>0$. Dually, an injective $A$-module $I$ is said to be \emph{$\nu^{-}$-stable} if $\nu_A^{-i}(I)$ are injective for all $i>0$.  The full subcategory of $A\pmodcat$ consisting of all $\nu$-stable projective $A$-modules is denoted by  $\mathscr{E}(A)$. Let $\epsilon_A$ be a basic additive generator for $\mathscr{E}(A)$, that is, $\epsilon_A$ is a basic module such that $\mathscr{E}(A)=\add(\epsilon_A)$. Recall that an module is basic if it is a direct sum of non-isomorphic indecomposable modules.

If $X$ is projective-injective and $\nu_A(X)\simeq X$ (or equivalently, $\nu_A^{-}(X)\simeq X$), then $X\in\mathscr{E}(A)$. So, the following lemma shows that $\epsilon_A$ is the maximal projective-injective basic $A$-module which generates $\mathscr{E}(A)$ and is closed under $\nu_A$ (or equivalently, under $\nu^{-}_A$).

\begin{Lem}\label{stp}
$(1)$ $\epsilon_A\in A\PI$ and $\nu_A(\epsilon_A)\simeq \epsilon_A\simeq \nu_A^{-}(\epsilon_A)$. In particular, the algebra $\End_A(\epsilon_A)$ is self-injective.

$(2)$ An $A$-module is $\nu$-stably projective if and only if it is $\nu^-$-stably injective.

$(3)$ The functor $\Hom_A(-,A)$ induces a duality
from $\mathscr{E}(A)$ to $\mathscr{E}(A^{\opp})$ which sends $\varepsilon_A$ to $\varepsilon_{A^{\opp}}$.
\end{Lem}

{\it Proof.}
$(1)$ Since  $\mathscr{E}(A)=\add(\epsilon_A)$, we see that $\nu_A^j(\epsilon_A)$
are projective for all $j\geq 0$. This implies that $\nu_A(\epsilon_A)\in\mathscr{E}(A)$. As $\epsilon_A$ is basic, the module $\nu_A(\epsilon_A)$ is isomorphic to a direct summand of
$\epsilon_A$. However, both $\nu_A(\epsilon_A)$ and $\epsilon_A$ have the same number of
indecomposable direct summands. Thus $\nu_A(\epsilon_A)\simeq \epsilon_A$. This leads to
$\epsilon_A\in A\PI$ and therefore $\epsilon_A\simeq \nu_A^{-}(\epsilon_A)$.
The last statement in $(1)$ is a result of Martinez-Villa (see, for example, \cite[Lemma 3.1(3)]{APT}).

$(2)$ Let $_AV$ be a basic additive generator for the category $\mathscr{V}$ of $\nu^-$-stably injective $A$-module. We can show that $_AV\in A\PI$ and $\nu_A^{-}(V)\simeq V$. This is dual to $(1)$. Further, for a projective-injective $A$-module $X$, it is known that $\nu_A(X)\simeq X$ if and only if $\nu_A^{-}(X)\simeq X$. Thus $\epsilon_A\in\mathscr{V}$ and $V\in\mathscr{E}(A)$. It follows that $\mathscr{E}(A)=\mathscr{V}$ and therefore $\epsilon_A\simeq V$.

$(3)$ This follows from $(1)$ and the definition of $\nu_A$.  $\square$

\medskip
Recall that the module $X$ is called a \emph{generator} over $A$ if $\add(_AA)\subseteq
\add(X)$; a \emph{cogenerator}  if $\add(D(A_A))\subseteq \add(X)$, and a
\emph{generator-cogenerator} if it is both a generator and a cogenerator
over $A$.

Let $V$ be a generator over $A$ with $B:=\End_A(V)$. Then $\Hom_A(V,I)$ is an injective $B$-module for every injective $A$-module $I$.
If $V$ is a generator-cogenerator, then each projective-injective $B$-module is precisely of the form $\Hom_A(V,I)$ with $I$ an injective $A$-module.
This is due to the isomorphism $D\Hom_A(P,V)\simeq \Hom_A(V,\nu_A(P))$ for all $P\in \add(_AA)$.

\medskip
The following observation may be useful to determine the dominant dimensions of modules.

\begin{Lem}
Let $X$ be an $A$-module with an exact sequence $0\ra X\ra X_0\ra \cdots\ra X_n\ra V\ra 0$ such that all $X_j$ are projective-injective modules. Suppose that $0\ra X\ra E_0\ra \cdots\ra E_n\lraf{d_n}E_{n+1}\ra \cdots$ is a minimal injective resolution of $X$. Then $E_{n+1}$ is projective if and only if so is the injective envelope of $V$. \label{inj-env}
\end{Lem}

{\it Proof.} Since $X_j$, with $0\le j\le n$, are projective-injective, we see from homological algebra that $E_j$, with $0\le j\le n$, are projective-injective. It follows from the dual version of Schanuel's Lemma that $V\oplus E_n\oplus X_{n-1}\oplus \cdots \oplus C \simeq \Img(d_n)\oplus X_n\oplus E_{n-1}\oplus \cdots\oplus C'$, where $C=X_0$ and $C'=E_0$ if $n$ is odd, and $C=E_0$ and $C'=X_0$ if $n$ is even. Thus, by taking injective envelopes, we obtain the following isomorphism of modules:
$$ I(V)\oplus Q\simeq E_{n+1}\oplus Q'$$
where $Q$ and $Q'$ are projective-injective modules. This implies that $E_{n+1}$ is projective if and only if $I(V)$ is projective. $\square$

\subsection{Derived equivalent algebras with different dominant dimensions}

In this section, we shall give a construction to produce derived equivalent algebras with different dominant dimensions.

As usual, for two $A$-modules $X$ and $Y$, we denote by $\mathscr{P}(X, Y)$ the set of homomorphisms from $X$ to $Y$ that factorise through a projective $A$-module, and by $\StHom_A(X, Y)$ the Hom-set in the stable category of $A\modcat$.
If $A$ is self-injective, then $\Ext^1_A(X, Y)\simeq
\StHom_A(\Omega_A(X), Y)$.

Our first main result in this paper is as follows:

\begin{Theo}\label{main result}
Let $A$ be a self-injective algebra over a field $k$. Suppose that
$$0\lra X\lra P\lra Y\lra 0$$ is an exact sequence of $A$-modules
with $P$ projective. Let $N$ be an $A$-module without non-zero projective
direct summands, and let
$$\small{\Lambda :=\left(
 \begin{array}{ccc}
\End_A(A)      & \Hom_A(A, N) & \Hom_A(A, X) \\
\Hom_A(N, A)   & \End_A(N)    & \Hom_A(N, X)\\
\Hom_A(X, A)   &\mathscr{P}(X, N) & \End_A(X)  \\
   \end{array}
 \right)\;\,\mbox{and}\;\,
\Gamma :=\left(
\begin{array}{ccc}
\End_A(A)      & \Hom_A(A, N)  & \Hom_A(A, Y) \\
\Hom_A(N, A)   & \End_A(N)     & \mathscr{P}(N, Y)\\
\Hom_A(Y, A)   & \Hom_A(Y, N)  & \End_A(Y)  \\
   \end{array}
 \right).}
$$
Then the following statements are true:

$(1)$ The algebras $\Lambda$ and $\Gamma$ are derived equivalent.

$(2)$ If $\StHom_A(N, Y)\neq 0$, then $\dm{(\Gamma)}=1$.

$(3)$ If $\StHom_A(X, N)\neq 0$, then $\dm{(\Lambda)}=1$.
\end{Theo}

{\it Proof.} $(1)$ In Proposition \ref{gder}, we take $\mathcal{C}:=
A\modcat$ and $M:= {}_AA\oplus N$. Since $P\in\add(_AA)\subseteq
\add(M)$, the given exact sequence
$$0\lra X\lra P\lra Y\lra 0$$ is an
$\add(M)$-exact sequence in $A\modcat$ (see Definition
\ref{ass}). Clearly, this sequence is always an $\add(_AA)$-split
sequence since $A$ is self-injective, but it does not have to be an
$\add(M)$-split sequence in general. So we cannot use Proposition \ref{der}. Now, one can check straightforward
that the rings $\Lambda$ and $\Gamma$ in Theorem \ref{main result}
are the same as the rings $R$ and $S$ in Proposition \ref{gder},
respectively. Thus $(1)$ follows from Proposition \ref{gder}.

$(2)$ Set $W:=A\oplus N\oplus Y$ and $B:=\End_A(W)$. In the sequel,
we always identify $B$ with the following matrix ring of $3\times 3$ matrices:
$$\left(
 \begin{array}{ccc}
\End_A(A)      & \Hom_A(A, N) & \Hom_A(A, Y) \\
\Hom_A(N, A)   & \End_A(N)    & \Hom_A(N, Y)\\
\Hom_A(Y, A)   & \Hom_A(Y, N) &  \End_A(Y)  \\
   \end{array}
 \right).
$$
Furthermore, let
 $$e_1:=\left(
 \begin{array}{ccc}
 1   & 0 & 0\\
 0   & 0 & 0 \\
 0   & 0 & 0 \\
   \end{array}
 \right), \;e_2:=\left(
 \begin{array}{ccc}
 0   & 0 & 0\\
 0   & 1 & 0 \\
 0   & 0 & 0 \\
   \end{array}
 \right)\; \mbox{and}\; e_3:=\left(
 \begin{array}{ccc}
 0   & 0 & 0\\
 0   & 0 & 0 \\
 0   & 0 & 1 \\
   \end{array}
 \right).
$$
Then
$Be_1\simeq \Hom_A(W, A)$ as $B$-modules. Moreover, the following statements hold:

$(a)$ The algebra $\Gamma$ is a subalgebra of $B$ and contains all
$e_i$ for $1\leq i\leq 3$.

$(b)$ $\Gamma e_1=  Be_1$ and $ \Gamma e_2= Be_2$.

$(c)$ There is an exact sequence
$$
0\lra \Gamma e_3 \lra Be_3\lra Be_3/\Gamma e_3\lra 0
$$
of $\Gamma$-modules such that $Be_3/\Gamma e_3\simeq \Hom_A(N,
Y)/\mathscr{P}(N, Y)=\StHom_A(N, Y)$ as $k$-modules,
$e_2(Be_3/\Gamma e_3)= Be_3/\Gamma e_3$ and $e_1(Be_3/\Gamma
e_3)=0=e_3(Be_3/\Gamma e_3)$.

First of all, we point out that $\nu_\Gamma(\Gamma e_1)\simeq
\nu_B(Be_1)$ as $\Gamma$-modules.

In fact, it follows from $e_1B=e_1\Gamma$ that there are the
following isomorphisms of $\Gamma$-modules:
$$\nu_\Gamma(\Gamma e_1)=D\Hom_\Gamma(\Gamma e_1, \Gamma)
\simeq D(e_1\Gamma)=D(e_1B)\simeq D\Hom_B(Be_1, B)\simeq
\nu_B(Be_1).$$

Next, we claim that $\add(_\Gamma \Gamma e_1)=\add\big(\nu_\Gamma(\Gamma
e_1)\big)\subseteq \Gamma\modcat$. This implies that $_\Gamma\Gamma e_1$
is projective-injective.

Actually, we always have
$\add(_BBe_1)=\add\big(\nu_B(Be_1)\big)\subseteq B\modcat$. To see this, we first note the following isomorphisms of
$B$-modules:
$$\nu_B(Be_1)=D\Hom_B(Be_1, B)\simeq D(e_1B)\simeq D(\Hom_A(A, W))\simeq \Hom_A(W, DA)$$
where the last isomorphism is due to the following well known
result: Let $C$ be a $k$-algebra and $P$ a $C$-module. If $P$ is
projective, then $D\Hom_C(P, U)\simeq \Hom_C(U, \nu_C(P))$ as
$k$-modules for each $C$-module $U$. Since the algebra $A$ is
self-injective, we certainly have $\add(_AA)=\add(_ADA).$
Particularly, this implies that $\add(_B\Hom_A(W,
A))=\add(_B\Hom_A(W, DA))$, and therefore
$\add(_BBe_1)=\add\big(\nu_B(Be_1)\big)\subseteq B\modcat$; and the module
$_BBe_1$ is projective-injective.

To prove $\add(_\Gamma \Gamma e_1)=\add\big(\nu_\Gamma(\Gamma
e_1)\big)$, we first observe the following fact:

Let $S\to R$ be a ring homomorphism of two $k$-algebras $S$ and $R$.
Suppose that $M$ and $L$ are $R$-modules. If $\add(_RM)=\add(_RL)$,
then $\add(_SM)=\add(_SL)$.

In our case, we know that $\Gamma\subseteq B$ is an extension of
algebras and that $\Gamma e_1=Be_1$ and $\nu_\Gamma(\Gamma
e_1)\simeq \nu_B(Be_1)$ as $\Gamma$-modules. By the above-mentioned
fact, we conclude from $\add(_BBe_1)=\add\big(\nu_B(Be_1)\big)\subseteq
B\modcat$ that $\add(_\Gamma \Gamma e_1)=\add\big(\nu_\Gamma(\Gamma
e_1)\big)\subseteq \Gamma\modcat$. This finishes the claim.

Now, we show $\dm(\Gamma)\geq 1$ .

Let $f: W\to I(W)$ be an injective envelope of $_AW$. Then
$I(W)\in\add(_ADA)=\add(_AA)$. Applying $\Hom_A(W,-)$ to the map
$f$, we obtain the induced map $f^*: B=\End_A(W)\to\Hom_A(W, I(W))$
which is injective. Clearly, $\Hom_A(W, I(W))\in\add(_B\Hom_A(W,
A))=\add(_BBe_1)$.

Recall that
$Be_1=\Gamma e_1$ which is a projective-injective $\Gamma$-module.
This implies that $\Hom_A(W, I(W))$ is also projective-injective as
a $\Gamma$-module. Since $\Gamma$ can be embedded into $B$ and
$f^*$ is an injective homomorphism of $\Gamma$-modules, there is
an injection $\Gamma \to \Hom_A(W, I(W))$ of $\Gamma$-modules. It
follows that $\dm(\Gamma)\geq 1$.

Finally, we assume $\StHom_A(N, Y)\neq 0$, and want to show
$\dm(\Gamma)=1$.

Let $g: Y\to I(Y)$ be an injective envelope of $_AY$. Then the map
$g^*: \Hom_A(W, Y)\to \Hom_A(W, I(Y))$ induced by $g$ is injective.
Since $Be_3\simeq \Hom_A(W, Y)$ as $B$-modules, there exists an
injection $\varphi: Be_3\to\Hom_A(W, I(Y))$ of $B$-modules.
Thus, using the sequence in $(c)$, we can construct the
following exact commutative diagram of $\Gamma$-modules:
$$
\xymatrix{ 0\ar[r] & \Gamma e_3\ar[r]\ar@{=}[d] &
Be_3\ar[r]\ar[d]^-{\varphi}
& Be_3/\Gamma e_3\ar[r]\ar@{-->}[d]^-{\psi} &0\\
0\ar[r] & \Gamma e_3\ar[r] & \Hom_A(W, I(Y))\ar[r] & V\ar[r] &0 }
$$
where $V$ is a $\Gamma$-module  and $\psi$ is injective. This means
that $Be_3/\Gamma e_3$ is isomorphic to a submodule of $_\Gamma V$,
and therefore $\soc(_\Gamma Be_3/\Gamma e_3)$ is a direct summand of
$\soc(_\Gamma V)$.

Since $\StHom_A(N, Y)\neq 0$ by assumption, we know from $(c)$ that
$Be_3/\Gamma e_3\neq 0$ and $e_1(Be_3/\Gamma e_3)=0=e_3(Be_3/\Gamma
e_3)$. It follows that there exists a simple $\Gamma$-module $S$,
which is a direct summand of $\soc(_\Gamma Be_3/\Gamma e_3)$, such
that its projective cover $P(S)$ belongs to $\add(_\Gamma \Gamma e_2)$.
Consequently, the module $S$ is a direct summand of $\soc(_\Gamma
V)$.

Since $I(Y)\in\add(_AI(W))$ and $_\Gamma\Hom_A(W, I(W))$ is
projective-injective, we see that $\Hom_A(W, I(Y))$ is
projective-injective as an $\Gamma$-module. So, by Lemma \ref{inj-env}, to prove
$\dm(\Gamma)=1$, it is sufficient to show that the injective
envelope $I(V)$ of the module $_\Gamma V$ is not projective. As $S$
is a direct summand of $\soc(_\Gamma V)$, we know that the injective
envelope $I(S)$ of $S$ is also a direct summand of $I(V)$. In the
following, we shall show that $I(S)$ is not projective.

Suppose that $_\Gamma I(S)$ is projective. Then $_\Gamma I(S)$ is
projective-injective. On the one hand, by the proof of
$\dm(\Gamma)\geq 1$, we have $ I(S)\in\add(_\Gamma \Gamma
e_1)=\add(\nu_\Gamma(\Gamma e_1)).$ It follows from
$I(S)=\nu_\Gamma(P(S))$ that $P(S)\in\add(_\Gamma \Gamma e_1)$. This
implies that $e_1P(S)\in\add(e_1\Gamma e_1)\subseteq e_1\Gamma
e_1\modcat$. On the other hand, since $P(S)\in \add(_\Gamma \Gamma
e_2)$, we obtain $e_1 P(S)\in\add(e_1\Gamma e_2)\subseteq e_1\Gamma
e_1\modcat$. Note that $e_1\Gamma e_1\simeq \End_A(A)\simeq A$ as
algebras, and that $e_1\Gamma e_2\simeq\Hom_A(A, N)\simeq N$ as
$A$-modules. Consequently, after identifying $e_1\Gamma e_1$ with
$A$, we see that $e_1P(S)\in\add(_AA)$ and $e_1P(S)\in \add(_AN)$.
This contradicts to the assumption that $_AN$ has no nonzero
projective modules as direct summands. Thus $_\Gamma I(S)$ is not
projective, and $\dm(\Gamma)=1$. This completes the proof of $(2)$.

$(3)$ Since the algebra $A$ is self-injective, this Hom-functor
$\Hom_A(-, A): A\modcat\to A\opp\modcat$ is a duality with dual
inverse $\Hom_{A\opp}(-, A): A\opp\modcat\to A\modcat$ (see
\cite[IV. Proposition 3.4]{ARS}). By the assumption on the module
$_AN$ in Theorem \ref{main result}, we infer that $\Hom_A(N, A)_A$
does not have any projective direct summands.

To prove $(3)$, we shall focus on the inclusion $\Lambda
\to\End_A(A\oplus N\oplus X)$ of algebras, and consider the right
$\Lambda$-modules instead of left $\Lambda$-modules. In this situation, we can show $\dm(\Lambda\opp)=1$ by following
the proof of $(2)$. Here, we leave the
details to the reader. Note that $\dm(\Lambda)=\dm(\Lambda\opp)$.
Thus $\dm(\Lambda)=1$. This finishes the proof of $(3)$. $\square$

\medskip
As a consequence of Theorem \ref{main result}, we have the following
result.
\begin{Koro}\label{cor1}
Let $A$ be a self-injective algebra. Suppose that $Y$ and $N$ are
$A$-modules such that $N$ has no non-zero projective direct
summands. If $\,\Ext^1_A(Y, N)=0$ and $\Ext^1_A(N, \Omega_A(Y))\neq
0$, then the endomorphism algebra $\End_A(A\oplus N\oplus
\Omega_A(Y))$ has  dominant dimension at least $2$ and is derived
equivalent the following matrix algebra of dominant dimension $1$
$$\left(
 \begin{array}{ccc}
\End_A(A)      & \Hom_A(A, N) & \Hom_A(A, Y) \\
\Hom_A(N, A)   & \End_A(N)    & \mathscr{P}(N, Y)\\
\Hom_A(Y, A)   & \Hom_A(Y, N) &  \End_A(Y)  \\
   \end{array}
 \right).$$
\end{Koro}

{\it Proof.} Clearly, the dominant dimension of $\End_A(A\oplus N\oplus
\Omega_A(Y))$ is at least $2$ (actually, it is equal to $2$ by the remarks after the definition of dominant dimensions). Since $A$ is self-injective, it is known that
$\StHom_A(\Omega_A(L), M)\simeq \Ext^1_A(L, M)$ and $\StHom_A(L,
M)\simeq \StHom_A(\Omega_A(L), \Omega_A(M))$ for $L, M\in A\modcat$.
Clearly, this implies that $\StHom_A(\Omega_A(Y), N)=0$ if and only
if $\Ext^1_A(Y, N)=0$, and that $\StHom_A(N, Y)\neq 0$ if and only
if $\Ext^1_A(N, \Omega_A(Y))\neq 0$.

Recall that there always exists an exact sequence of $A$-modules:
$$0\lra \Omega_A(Y)
\lra P\lraf{\pi} Y\lra 0$$  such that $\pi$ is a projective cover of
$_AY$. Now, we take $X:=\Omega_A(Y)$ in Theorem \ref{main result} and get Corollary \ref{cor1} immediately from Theorem
\ref{main result}. $\square$

\medskip
Algebras of the form $\End_A(A\oplus M)$ with $A$ a self-injective algebra and $M\in A\modcat$ are called \emph{Morita algebras} in \cite{KY}. The above corollary shows that both Morita algebras and dominant dimensions are not invariant under derived equivalences, though self-injective algebras are invariant under derived equivalences.

We observe the following characterization of Morita algebras, its proof follows directly from \cite[Corollary 1.4, Theorem 1.5]{KY}.

\begin{Lem}\label{Morita}
Let $A$ be an algebra of dominant dimension at least $2$. Then $A$ is a Morita algebra if and only if $A\PI\, =\mathscr{E}(A)$.
\end{Lem}

%

\section{Tilting modules, dominant dimensions and Morita algebras \label{sect4new}}

Tilting modules supply an important class of derived equivalences.
In the following we consider when derived equivalences given by tilting modules preserve dominant dimensions and Morita algebras. During the course of our discussions, we also establish a lower bound for dominant dimensions of algebras under tilting procedure.

\subsection{Dominant dimensions for general tilting modules }
\medskip

Let us first recall the definition of tilting modules (see, for instance, \cite{H} or any text book on the representation theory of finite dimensional algebras).

An $A$-module $T$ is called a \emph{tilting module} if $\pd(_AT)=n<\infty$, $\Ext^j_A(T,T)=0$ for all $j>0$, and there is an exact sequence
$ 0\ra {}_AA\ra X_0\ra X_1\ra \cdots \ra X_n\ra 0$ in $A$-mod with all $X_j\in \add(T)$.

Clearly, by definition, indecomposable projective-injective modules are isomorphic to direct summands of each tilting module.

\medskip
From now on, we investigate behaviours of dominant dimensions of the endomorphism algebras of tilting modules.

Let $T$ be a tilting $A$-module of projective dimension $n\geq 1$, and let
$B:=\End_A(T)$. We first fix a minimal projective resolution of $_AT$ as follows:
$$
0\lra P_n\lra P_{n-1}\lra\cdots\lra P_1\lra P_0\lra T\lra 0.
$$
Then $A\pmodcat\,=\add(\bigoplus_{i=0}^{n}P_i)$, and any projective summand of $T$ belongs to $\add(P_0)$. Further, $T_B$ is a tilting right $B$-module with $\pd(T_B) = n$. It is well known that $A$ and $B$ are derived equivalent (see \cite{Happel}).

For convenience, we introduce the following definition which seems to be useful in the rest of our discussions.

\begin{Def}\label{heart}
Let $_AT=P\oplus T'$ such that $P$ is projective and $T'$ has no non-zero projective direct summands. The \emph{heart} of $T$ is defined to be a basic $A$-module $E(A,T)$ such that
$$\add\big(E(A,T)\big)=\{X\in\add(_AP)\mid \nu_A(X)\in\add(_AT)\}.$$
\end{Def}

Note that if $T={_A}A$, then $\nu_A(E(A,T))$ is an additive generator for $A\PI$. In general, the module $E(A,T)$ may not be injective. Since $\nu_A(\epsilon_A)\simeq \epsilon_A \in A\PI \, \subseteq \add(_AT)$ by Lemma \ref{stp}, we always have $\mathscr{E}(A)\subseteq\add(E(A,T))$. Moreover, if $_AT'$ has no non-zero injective direct summands, then $\add(E(A,T))=\{X\in\add(_AP)\mid\nu_A(X)\in\add(_AP)\}.$

\smallskip
Throughout this section, we fix a tilting $A$-module $T$ with the above projective resolution, and let $P$ and $E(A,T)$ be defined as above. If no confusion arises, we simply write $E$ for $E(A,T)$.

\medskip
The following homological fact will be used for later discussions.
\begin{Lem}\label{iso}
Let $M$ and $N$ be $A$-modules. Then we have the following:

$(1)$ If $N\in\add(_AA)$, then the functor $\Hom_A(-,T)$ induces an isomorphism of abelian groups: $\Hom_A(M,N)\simeq \Hom_{B^{\opp}}(\Hom_A(N,T),\Hom_A(M,T))$.

$(2)$ If $M\in \add(D(A_A))$, then the functor $\Hom_A(T,-)$ induces an isomorphism of abelian groups: $\Hom_A(M,N)\simeq \Hom_{B}(\Hom_A(T,M),\Hom_A(T,N))$.
\end{Lem}

{\it Proof.}  To show $(1)$, we use the fact that $_AA$ has an $\add(_AT)$-copresentation, that is, there is an exact sequence $0\to A\to T_0\to T_1$ of $A$-modules with $T_0,T_1\in\add(_AT)$ such that the sequence $\Hom_A(T_1,T)\to \Hom_A(T_0,T)\to \Hom_A(A,T)\to 0$ is still exact. Dually, to show $(2)$, we use the fact that $D(A_A)$ has an $\add(_AT)$-presentation, that is, there is an exact sequence $T_1'\to T_0'\to D(A_A)\to 0$ of $A$-modules with $T_0',T_1'\in\add(_AT)$ such that the sequence $\Hom_A(T,T_1')\to \Hom_A(T,T_0')\to \Hom_A(T,D(A))\to 0$ is still exact. $\square$

\medskip
The following result gives a characterization of projective-injective $B$-modules.

\begin{Lem} \label{proj-inj}
The functor $\Hom_A(-,T):A\modcat\to B^{\opp}\modcat$ restricts to the following two dualities between additive categories:
$$\add(E)\lraf{\simeq} B^{\opp}\PI\;\;\mbox{and}\;\;\add(\mathscr{E}(A))\lraf{\simeq}
\add(\mathscr{E}(B^{\opp}))$$
which send $\varepsilon_A$ to $\varepsilon_{B^{\opp}}$. In particular, we have $\End_A(\varepsilon_A)\simeq \End_{B}(\varepsilon_{B})$ as algebras.
\end{Lem}

{\it Proof.} For simplicity, we set $F:=\Hom_A(-,T)$.  Clearly, for any $X\in\add(E)$, we have $F(X)\in B^{\opp}\pmodcat$ due to $X\in\add(P)\subseteq \add(_AT)$. Since $X$ is projective and $\nu_A(X)\in\add(_AT)$, we see that $$DF(X)=D\Hom_A(X,T)\simeq \Hom_A(T,\nu_A(X))\in B\pmodcat\, .$$  This forces $F(X)\in B^{\opp}\imodcat$, and therefore  $F(X)\in B^{\opp}\PI$. So, $F$ induces a functor $F_1:\add(E)\to B^{\opp}\PI$. Since $\add(E)\subseteq\add(_AA)$, the functor $F_1$ is fully faithful by Lemma \ref{iso}(1).

Now, we show that if $Q$ is an indecomposable projective-injective $B^{\opp}$-module, then there exists an $A$-module $Y\in\add(E)$ such that $Q\simeq F(Y)$. This will verify that $F_1$ is dense, and therefore $F_1$ is a duality.

In fact, from the surjective map $P_0\to T$ we obtain an injective homomorphism $B\to F(P_0)$ of $B^{\opp}$-modules. It follows that $Q$ is isomorphic to a direct summand of $F(P_0)$.
Note that $F$ induces a ring isomorphism from $\End_A(P_0)^{\opp}$ to
$\End_{B^{\opp}}(F(P_0))$ by Lemma \ref{iso}(1) due to $P_0\in\add(A)$. Consequently, there is a direct summand $Y$ of $P_0$ such that $Q_B\simeq F(Y)$. It suffices to check
that $Y\in\add(_AT)$ and $\nu_A(Y)\in\add(_AT)$ since $\add(P)=\add(P_0)\cap\add(T)$.

On the one hand, since $F(Y)_B$ is projective, there exists a module $T'\in\add(_AT)$ such that $F(Y)\simeq F(T')$. We claim that $Y\simeq T'$. Actually, since $T'\in\add(_AT)$, we first have $\Hom_A(Y,T')\simeq\Hom_{B^{\opp}}(F(T'),F(Y))$. Further, since $Y\in\add(P_0)$, we then see from Lemma \ref{iso}(1) that $\Hom_A(T',Y)\simeq\Hom_{B^{\opp}}(F(Y),F(T'))$. Thus $Y\simeq T'\in\add(_AT)$.

On the other hand, since $_AY$ is projective and $F(Y)_B$ is injective, we have $$DF(Y)=D\Hom_A(Y,T)\simeq \Hom_A(T,\nu_A(Y))\in \pmodcat B .$$
There exists a module $T''\in\add(_AT)$ such that $G(\nu_A(Y))\simeq G(T'')$, where $G$ is the functor $\Hom_A(T,-):A\modcat\to B\modcat$. We claim that $\nu_A(Y)\simeq T''$. In fact, since $T''\in\add(_AT)$, we have $\Hom_A(T'',\nu_A(Y))\simeq\Hom_{B}\big(G(T''),G(\nu_A(Y))\big)$. As $\nu_A(Y)\in \add(D(A_A))$, we see from Lemma \ref{iso}(2) that $\Hom_A(\nu_A(Y),T'')\simeq\Hom_{B}(G(\nu_A(Y)),G(T''))$. Therefore $\nu_A(Y)\simeq T''\in\add(_AT)$. Thus we have shown that $F_1$ is a duality.

Next, we show that $F_1$ induces a duality $F_2: \add(\mathscr{E}(A))\to\add(\mathscr{E}(B^{\opp}))$.

Recall that $\add(\mathscr{E}(A))=\add(\varepsilon_A)\subseteq A\PI$ and $\add(\mathscr{E}(B^{\opp}))=\add(\varepsilon_{B^{\opp}})\subseteq \PI {B^{\opp}}$. Since $\nu_A(\varepsilon_A)\simeq \varepsilon_A$ by Lemma \ref{stp}(1), we have $\add(\mathscr{E}(A))\subseteq\add(E)$. Note that, for a projective $A$-module $Z$, the following isomorphisms are true: $$\nu_{B^{\opp}}^{-}(F(Z))=\Hom_B(DF(Z),B)\simeq \Hom_B\big(G(\nu_A(Z)),B\big)\simeq F(\nu_A(Z)).$$
It follows that $\nu_{B^{\opp}}^{-}(F(Z))\simeq F(Z)$ if and only if $ F(\nu_A(Z))\simeq F(Z)$.
In particular, if $\nu_A(Z)\simeq Z$, then $\nu_{B^{\opp}}^{-}(F(Z))\simeq F(Z)$. Since $\varepsilon_A\in\add(_AA)$ and $\nu_A(\varepsilon_A)\simeq \varepsilon_A$, we have $\nu_{B^{\opp}}^{-}(F(\varepsilon_A))\simeq F(\varepsilon_A)$. This implies that $F(\varepsilon_A)\in\add(\mathscr{E}(B^{\opp}))$ by Lemma \ref{stp}(1), and therefore $F(\varepsilon_A)$ is isomorphic to a direct summand of $\varepsilon_{B^{\opp}}$. So, the functor $F_1$ restricts to a fully faithful functor $F_2:\add(\mathscr{E}(A))\to\add(\mathscr{E}(B^{\opp}))$.

It remains to show that $F_2$ is dense. Indeed, since $F_1$ is a duality and since $\varepsilon_{B^{\opp}}\in B^{\opp}\PI$ is basic, there exists a basic $A$-module $E'\in\add(E)$ such that $F(E')\simeq \varepsilon_{B^{\opp}}$ as $B^{\opp}$-modules.
As $\nu_{B^{\opp}}^{-}(\varepsilon_{B^{\opp}})\simeq \varepsilon_{B^{\opp}}$ by Lemma \ref{stp}(1), we have $\nu_{B^{\opp}}^{-}(F(E'))\simeq F(E')$. Thus $F(\nu_A(E'))\simeq F(E')$.
We claim that $\nu_A(E')\simeq E'$. Actually, since $E'\in\add(E)$, we have $\nu_A(E')\in\add(_AT)$. This gives rise to  $\Hom_A(E',\nu_A(E'))\simeq\Hom_{B^{\opp}}(F(\nu_A(E')),F(E'))$.
Further, since $E'\in\add(E)\subseteq\add(_AA)$, we see from Lemma \ref{iso}(1) that $\Hom_A(\nu_A(E'),E')\simeq\Hom_{B^{\opp}}\big(F(E'),F(\nu_A(E'))\big)$. Thus $\nu_A(E')\simeq E'$. By Lemma \ref{stp}(1), the module $E'$ is isomorphic to a direct summand of $\varepsilon_A$,
and therefore $\varepsilon_{B^{\opp}}$ is isomorphic to a direct summand of $F(\varepsilon_A)$. Hence, $F(\varepsilon_A)\simeq\varepsilon_{B^{\opp}}$ and $F_2$ is a duality. $\square$

\medskip
The following corollary shows when the endomorphism algebra of a tilting module is a Morita algebra.
\begin{Koro}\label{B-M}
If $\dm(B)\ge 2$, then $B$ is a Morita algebra if and only if
$\mathscr{E}(A)=\add(E)$.
\end{Koro}

{\it Proof.} By Lemma \ref{Morita}, we see that $B$ is a Morita algebra if and only if $B\PI\, =\add(\varepsilon_B)$. This is also equivalent to that $B^{\opp}\PI\, =\add(\varepsilon_{B^{\opp}})$ by Lemma \ref{stp}(1) and (3).
Now, Corollary \ref{B-M} follows from Lemma \ref{proj-inj}. $\square$

\medskip
As an consequence of Corollary \ref{B-M}, we get a class of tilting modules which transfer Morita algebras again to Morita algebras.

\begin{Koro}\label{app}
Suppose that $A$ is a Morita algebra and $\dm(B)\ge 2$.

$(1)$ If the non-projective part $_AT'$ of $T$ has no injective direct summands, then $B$ is a Morita algebra.

$(2)$ If $\dm(A)\geq n+1$, then $B$ is a Morita algebra. In particular, if $\pd(_AT)=1$, then $B$ is a Morita algebra.
\end{Koro}

{\it Proof.} $(1)$ Clearly, we have $\mathscr{E}(A)\subseteq\add(E)$. Suppose that $_AT'$ has no injective direct summands. Then $\add(E)=\{X\in\add(_AP)\mid\nu_A(X)\in\add(_AP)\}$. In particular, we obtain $\nu_A(E)\in\add(_AP)\subseteq \add(_AA)$. Since $\nu_A(E)$ is injective, we have
$\nu_A(E)\in A\PI$. Note that $A$ is a Morita algebra by assumption. It follows from Lemma \ref{Morita} that $A\PI\,=\mathscr{E}(A)$. This implies $\nu_A(E)\in\mathscr{E}(A)$, and therefore $E\in\mathscr{E}(A)$ by Lemma \ref{stp}(1).
Thus $\add(E)=\mathscr{E}(A)$. Now $(1)$ is a consequence of Corollary \ref{B-M}.

$(2)$ Suppose $\dm(A)\geq n+1$. Since $\dm(A)=\dm(A^{\opp})$, each projective $A^{\opp}$-module $U$ has dominant dimension at least $n+1$. This implies that the injective dimension of $U$ is either $0$ or at least $n+1$. Dually, the projective dimension of each injective $A$-module is either $0$ or at least $n+1$. Moreover, since $_AT=P\oplus T'$ and $\pd(_AT)\le n$, the module
$_AT'$ has no injective direct summands. Now, $(2)$ follows from $(1)$. $\square$

\medskip
From now on, let $\omega$ be a basic $A$-module such that $A\PI\,=\add(\omega)$. Note that projective-injective $A$-modules always appear in the projective summands of each tilting $A$-module. So $\omega\in\add(_AP)$.

\begin{Lem}\label{dm-formula}
The following statements are true:

$(1)$ $\nu_A(E)\emph{-}\dm(A)\leq \dm(B)+n$.

$(2)$ $\Hom_A(\nu_A^{-}(\omega),T)\emph{-}\dm(B^{\opp})\leq \dm(A)+n$.
\end{Lem}

{\it Proof.} $(1)$ Since $E$ is projective, the $A$-module $\nu_A(E)$ is injective and $\nu_A(E)\emph{-}\dm(A)$ makes sense. Let $s:=\nu_A(E)\emph{-}\dm(A)$ and $t:=\dm(B)$. Clearly, if $s\leq n$ or $t=\infty$, then $(1)$ holds automatically. It remains to show that if $s\geq n+1$ and $t<\infty$, then $s\leq t+n$.

Suppose that $s\geq t+n+1$. Then $\nu_A(E)\emph{-}\dm(X)\geq s\geq t+n+1$ for any projective $A$-module $X$. So, for the  minimal projective resolution of ${_A}T$
$$
0\lra P_n\lra P_{n-1}\lra\cdots\lra P_1\lra P_0\lra T\lra 0,
$$
we have $\nu_A(E)\emph{-}\dm(P_i)\geq t+n+1$ for all $0\leq i\leq n$.
It follows from Corollary \ref{calculating dd}(2) that
$$\nu_A(E)\emph{-}\dm({_A}T)\geq\min\{\nu_A(E)\emph{-}\dm(P_i)+n-i-1\mid 0\leq i\leq n\}-n+1
\geq t+1.$$
Thus there exists an exact sequence of $A$-modules:
$$0\lra {_A}T\lra E^0\lra E^1\lra \cdots\lra E^{t-1}\lra E^t $$
such that $E^i\in\add(\nu_A(E))$ for all $0\leq i\leq t$. Note that $\Ext_A^j(T,T)=0$ for each $j\geq 1$ since ${_A}T$ is a tilting module. Applying $\Hom_A(T,-)$ to the above exact sequence, we obtain the following exact sequence of $B$-modules:
$$
0\lra {_B}B\lra \Hom_A(T, E^0)\lra \Hom_A(T, E^1)\lra\cdots \lra\Hom_A(T, E^{t-1})\lra
\Hom_A(T, E^t)
$$
such that $\Hom_A(T, E^i)\in\add\big(\Hom_A(T, \nu_A(E))\big)$ for all $0\leq i\leq t$.


Now, we point out that
$\add\big(\Hom_A(T, \nu_A(E))\big)=B\PI$.
In fact, by Lemma \ref{proj-inj}, the functor $\Hom_A(-, T)$ induces an equivalence from $\add(E)$ to $B^{\opp}\PI$. Note that $D\Hom_A(E,T)\simeq\Hom_A(T,\nu_A(E))$ as $B$-modules since ${_A}E$ is projective. Thus $\Hom_A(T,-)$ induces an equivalence from $\add(\nu_A(E))$ to $B\PI$.
From this, we conclude that $\dm(B)\geq t+1$, which is a contradiction. Thus $s\leq t+n$ and $(1)$ holds.

$(2)$ Since ${_A}\omega$ is injective, we see that $\nu_A^{-}(\omega)$ is projective.
It follows that $D\Hom_A(\nu_A^{-}(\omega),T)\simeq\Hom_A(T, \omega)$. Since $\omega\in\add(P)\subseteq\add(_AT)$,  the $B$-module $\Hom_A(T, \omega)$ is projective. So $\Hom_A(\nu_A^{-}(\omega),T)_B$ is injective and $\Hom_A(\nu_A^{-}(\omega),T)\emph{-}\dm(B^{\opp})$ makes sense.

Since ${_A}T$ is an $n$-tilting module and $B=\End_A(T)$, it is known that $T_B$ is also an $n$-tilting module and that $\End_{B^{\opp}}(T)\simeq A^{\opp}$. Let $E':=E(B^{\opp},T)$ be
the heart of $T_B$. We claim that $E'\simeq \Hom_A(\omega, T)$ as $B^{\opp}$-modules. Actually, by Lemma \ref{proj-inj}, the functor $\Hom_{B^{\opp}}(-, T)$ induces an equivalence from $\add(E')$ to $A\PI$. Since $A\PI\,=\add(\omega)$ and ${_A}\omega$ is basic, we have $\Hom_{B^{\opp}}(E', T)\simeq \omega$. It follows that $\Hom_A(\Hom_{B^{\opp}}(E', T), T)\simeq \Hom_A(\omega,T)$. Observe that $E'\simeq\Hom_A(\Hom_{B^{\opp}}(E', T), T)$ as $B^{\opp}$-modules since $E'$ is projective. Thus $E'\simeq \Hom_A(\omega, T)$. This verifies the claim. Since ${_A}\omega$ is projective-injective, we have
$$\nu_{B^{\opp}}(E')\simeq D\Hom_B(T,\omega)\simeq \Hom_A(\nu_A^{-}(\omega),T).$$
Now, by $(1)$, we see that $\Hom_A(\nu_A^{-}(\omega),T)\emph{-}\dm(B^{\opp})\leq \dm(A^{\opp})+n$. Since $\dm(A)=\dm(A^{\opp})$, the statement $(2)$ holds true. $\square$

\medskip
The following is our second main result which shows how the dominant dimensions of $A$ and $B$ are related.

\begin{Theo}\label{distance}
$(1)$ If $\omega\in \add(\nu_A(E))$, then $\dm(A)\leq \dm(B)+n$.

$(2)$ If $\nu_A(E)\in\add(\omega)$, then $\dm(B)\leq\dm(A)+n$.
\end{Theo}

{\it Proof.} Recall that $A\PI\,=\add(\omega)$. If $\omega\in \add(\nu_A(E))$, then
$A\PI\,\subseteq \add(\nu_A(E))$ and therefore $\dm(A)=\omega\emph{-}\dm(A)\leq \nu_A(E)\emph{-}\dm(A)$. By Lemma \ref{dm-formula}(1), we have
$\nu_A(E)\emph{-}\dm(A)\leq \dm(B)+n$. It follows that $\dm(A)\leq \dm(B)+n$. Thus $(1)$ holds.

Note that $B^{\opp}\PI\,=\add(\Hom_A(E,T))$ by Lemma \ref{proj-inj}. If
$\nu_A(E)\in\add(\omega)$, then $E\in \add(\nu_A^{-}(\omega))$ and therefore
$\dm(B^{\opp})= \Hom_A(E,T)\emph{-}\dm(B^{\opp})\leq
\Hom_A(\nu_A^{-}(\omega),T)\emph{-}\dm(B^{\opp}).$
Further, by Lemma \ref{dm-formula}(2), we see that  $\Hom_A(\nu_A^{-}(\omega),T)\emph{-}\dm(B^{\opp})\leq \dm(A)+n$. Since $\dm(B)=\dm(B^{\opp})$, we obtain $\dm(B)\leq \dm(A)+n$. Thus $(2)$ holds. $\square$

\medskip
From Theorem \ref{distance} we know that if $\add(\omega)= \add(\nu_A(E))$, then $\dm(A)=\infty$ if and only if $\dm(B)=\infty$.
As another consequence of Theorem \ref{distance}, we have the following result.

\begin{Koro}\label{compare-distance}
$(1)$ If $\nu_A(\omega)\simeq \omega$, then  $\dm(A)\leq \dm(B)+n$.

$(2)$ If the non-projective part $_AT'$ of $T$ has no injective direct summands, then $\dm(B)\leq \dm(A)+n$.

$(3)$ If $\add(\omega)=\add(\nu_A(E))$ (for example, $\omega\simeq \nu_A(E))$, then $|\dm(A)-\dm(B)|\leq n$.
\end{Koro}

{\it Proof.} If $\nu_A(\omega)\simeq \omega$, then $\omega\in\add(E)$ since $\omega\in\add(P)$. In this case, we have $\omega\in\add(\nu_A(E))$. Now, the statement $(1)$ follows from Theorem \ref{distance}(1).

If $_AT'$ has no non-zero injective direct summands, then
$$\add(E)=\{X\in\add(P)\mid \nu_A(X)\in\add(P)\}=\{X\in\add(P)\mid \nu_A(X)\in\add(\omega)\}$$
and therefore $\nu_A(E)\in \add(\omega)$. Thus $(2)$ follows from Theorem \ref{distance}(2).
Clearly, the statement $(3)$ is also due to Theorem \ref{distance}. $\square$

\begin{Koro}\label{compare-Morita}
If $A$ is a Morita algebra, then $\dm(A)\leq \dm(B)+n.$ In particular, if both $A$ and $B$ are Morita algebras, then $|\dm(A)-\dm(B)|\leq n$.
\end{Koro}

{\it Proof.} Since $A$ is a Morita algebra, we see from Lemma \ref{Morita} that $A\PI\,=\mathscr{E}(A)$. Note that $A\PI\,=\add(\omega)$, $\mathscr{E}(A)=\add(\epsilon_A)$ and $\nu_A(\epsilon_A)\simeq\varepsilon_A$ by Lemma \ref{stp}(1). Hence $\nu_A(\omega)\simeq\omega$. So the first statement of Corollary \ref{compare-Morita} follows from Corollary \ref{compare-distance}(1).

In addition, if $B$ is a Morita algebra, then we may use  the tilting right $B$-module $T_B$ of $\pd(T_B)= n$. In this case, we have $\dm(B^{\opp})\le \dm(\End_{B^{\opp}}(T_B))+n$ by the foregoing proof. Due to $\End_{B^{\opp}}(T_B)\cong A^{\opp}$, we obtain $\dm(B)=\dm(B^{\opp})\le \dm(A^{\opp})+n = \dm(A)+n$. Thus
$|\dm(A)-\dm(B)|\leq n$. $\square$

\medskip
Remarks. (1) If $A$ is a Morita algebra with $\dm(A)\geq n+2$, then $B$ is a Morita algebra. Indeed, it follows from Corollary \ref{compare-Morita} that $n+2\le \dm(A)\le \dm(B)+n$ and therefore $\dm(B)\ge 2$. Now, by Corollary \ref{app}(2), $B$ is a Morita algebra.

(2) If both $A$ and $B$ are assumed to be Morita algebras, then there is a general statement in \cite{FH}: For any derived equivalence between $A$ and $B$, there holds $\dm(B)-n\le \dm(A)\leq \dm(B)+n,$ where $n+1$ is the number of non-zero terms of a tilting complex defining the derived equivalence between $A$ and $B$. Unfortunately, our example in Section \ref{sect4} shows that $B$ in Corollary \ref{compare-Morita} does not have to be a Morita algebra.

\medskip
To give an optimal lower bound for the dominant dimension of $B$, we introduce the following definition of $T$-gradients of projective modules.

\begin{Def}
For a projective $A$-module $X$, let $$\cpx{T}_X:\quad \cdots\lra T^{-i}_X\lra T^{-i+1}_X\lra \cdots\lra T^{-2}_X\lra T^{-1}_X\lra \nu_A(X)\lra 0$$ be a minimal right $\add(_AT)$-approximation sequence of $\nu_A(X)$, where $T_X^0:=\nu_A(X)$.

$(1)$ The $T$-gradient of $X$, denoted by  $\partial_T(X)$, is defined as follows:
$$
\partial_T(X):=\inf\big\{i\geq 0\mid  T^{-i-1}_X\notin\add(\nu_A(E))\big\}.
$$

$(2)$ The $T$-gradient of the algebra $A$ is the $T$-gradient of $_AA$. Further, the global $T$-gradient of $A$ is
$$\partial(A,T):=\min\{\partial_T(P_i)+i\mid 0\leq i\leq n\}.$$

$(3)$ The tilting gradient of $A$ is $$\partial(A):= \sup\big\{\partial(A,T)\mid T \mbox{\; is a tilting \;} A\mbox{-module} \big\}.$$
\end{Def}

Remark that the sequence $\cpx{T}_X$ is always exact since an $A$-module $N$ with $\Ext_A^i(T,N)=0$ for all $i\geq 1$ is a quotient module of some module in $\add({_A}T)$. Thus one can define the $T$-gradients for all such modules $N$. However, in this paper, we confine our discussions to injective modules. Clearly, if $X\in\add(E)$, then $\partial_T(X)=\infty$.
Moreover, if $X, Y\in\add(_AA)$, then $\partial_T(X\oplus Y)=\min\{\partial_T(X), \partial_T(Y)\}$. Thus $\partial_T(A)\leq\partial_T(X)$ for all $X\in \add(_AA)$, and $\partial_T(A)\le \partial(A,T)\le \partial(A)$.

Note that $\dm(A)=\partial_{T}(A)=\partial(A,T)$ for $T={}_AA$ (see also Corollary \ref{equality}) since $\dm(A)=\dm(A^{\opp})$ and $\nu_A(E(A,A))$ is an additive generator for $A\PI$. In this sense, the definition of gradients generalizes the one of dominant dimensions of algebras.

Let us start with the following description of $T$-gradients.

\begin{Lem}\label{approx}
For any $X\in\add(_AA)$,  let  $$\cpx{E}_X:\quad \cdots\lra E^{-i}_X\lra E^{-i+1}_X\lra \cdots\lra E^{-2}_X\lra E^{-1}_X\lra X\lra 0$$ be a minimal right $\add(_AE)$-approximation sequence of $_AX$, where $E^0_X:=X$. Then

$$
\partial_T(X)=\dm(\Hom_A(X,T)_B)=\inf\{i\geq 0\mid \Hom_{\K A}\big(\cpx{E}_X, T[i]\big)\neq 0\}.
$$
\end{Lem}

{\it Proof.}  Recall that  ${_A}E\in\add(P)\subseteq\add(_AA)\cap\add(_AT)$ and that $B^{\opp}\PI\, =\add(\Hom_A(E,T))$ due to Lemma \ref{proj-inj}.

Let $s:=\dm(\Hom_A(X,T)_B)$ and $t:=\inf\{i\geq 0\mid \Hom_{\K A}\big(\cpx{E}_X, T[i]\big)\neq 0\}.$ Then, to show $s=t$, it suffices to prove that, for any integer $m\ge 1$, we have $s\ge m$ if and only if $t\ge m$.

Suppose that $s\geq m\geq 1$. Then $\Hom_A(X,T)_B$ has  a minimal injective resolution, starting with the following exact sequence
$$
(\dag)\quad 0\lra \Hom_A(X,T)\lra \Hom_A(E_1,T)\lra \Hom_A(E_2,T)\lra\cdots\lra\Hom_A(E_{m-1},T)\lra \Hom_A(E_{m}, T)
$$
such that $E_i\in\add(_AE)$ for $1\leq i\leq m$. Note that $\Hom_{B^{\opp}}(\Hom_A(U,T),\Hom_A(V,T))\simeq \Hom_A(V,U)$ for any projective $A$-modules $U$ and $V$ by Lemma \ref{iso}(1). So there exists a sequence
$$
\cpx{E}:\quad E_{m}\lraf{f_m} E_{m-1}\lraf{f_{m-1}} \cdots\lra E_2\lraf{f_2} E_1\lraf{f_1} X \lraf{f_0} 0
$$
such that $\cpx{\Hom}_A(\cpx{E},T)$ is precisely the sequence $(\dag)$. Since $_AT$ is a tilting module, we have an injection $A\to T_0$ for some $T_0\in\add(_AT)$. As both $E_i$ and $X$ belong to $\add(_AA)$, the exactness of the sequence $(\dag)$ implies that $f_{i+1}f_i=0$. Moreover, since $\nu_A(E)\in\add(_AT)$, the sequence $\cpx{\Hom}_A(\cpx{E},\nu_A(E))$ is exact. Thus it follows from $\cpx{\Hom}_A(\cpx{E},\nu_A(E))\simeq D\cpx{\Hom}_A(E,\cpx{E})$ that the sequence $\cpx{\Hom}_A(E,\cpx{E})$ is also exact. In other words, the homomorphism $g_i: E_i\to \Ker(f_{i-1})$  induced from $f_i$ is a right $\add(E)$-approximation of $\Ker(f_{i-1})$. Note that $g_i$ is minimal because the sequence $(\dag)$ is a part of a minimal injective resolution of $\Hom_A(X,T)_B$. Hence, $\cpx{E}$ can be regarded as a starting part of the minimal right $\add(E)$-approximation $\cpx{E}_X$ of $X$. So we may write the complex $\cpx{E}_X$ as follows:
$$\cpx{E}_X: \quad \cdots \lra E_X^{-m-2}\lra E_X^{-m-1}\lra  E_{m}\lraf{f_m} E_{m-1}\lraf{f_{m-1}} \cdots\lra E_2\lraf{f_2} E_1\lraf{f_1} X \lraf{f_0} 0,$$
and get $\Hom_{\K A}(\cpx{E}_X,T[i])\simeq H^i\big(\cpx{\Hom}_A(\cpx{E}_X, T)\big) =0$ for $0\leq i\leq m-1$. This implies $t\ge m$

Conversely, if $t\geq m\geq 1$, then $H^i\big(\cpx{\Hom}_A(\cpx{E}_X, T)\big)\simeq\Hom_{\K A}(\cpx{E}_X,T[i])=0$ for $0\leq i\leq m-1$, and therefore the following sequence
$$0\to \Hom_A(X,T)\lra \Hom_A(E_X^{-1},T)\lra \Hom_A(E_X^{-2},T)\lra\cdots\lra\Hom_A(E_X^{1-m},T)\lra \Hom_A(E_X^{-m}, T)$$
induced from $\cpx{E}_X$ is exact. Since $E_X^{-i}\in\add(E)$ for $1\leq i\leq m$, we have
$\dm(\Hom_A(X,T))\geq m$, that is, $s\ge m$. Thus $s=t$.

To check $\partial_T(X)=s$, we note from $\add(E)\subseteq\add(_AA)$ that $\cpx{\Hom}_A(\cpx{E}_X, T)\simeq D\cpx{\Hom}_A(T, \nu_A(\cpx{E}_X))$ as complexes.
It then follows that
$$\Hom_{\K A}(\cpx{E}_X,T[i])\simeq H^i\big(\cpx{\Hom}_A(\cpx{E}_X, T)\big)\simeq H^i\big(D\cpx{\Hom}_A(T, \nu_A(\cpx{E}_X))\big)\simeq DH^{-i}\big(\cpx{\Hom}_A(T, \nu_A(\cpx{E}_X))\big).$$ Hence $\Hom_{\K A}(\cpx{E}_X,T[i])=0$ if and only if
$H^{-i}\big(\cpx{\Hom}_A(T, \nu_A(\cpx{E}_X))\big)=0$, where the complex
$$\nu_A(\cpx{E}_X):\quad \cdots\lra \nu_A(E^{-i}_X)\lra \nu_A(E^{-i+1}_X)\lra \cdots\lra \nu_A(E^{-2}_X)\lra \nu_A(E^{-1}_X)\lra \nu_A(X)\lra 0$$
has terms $\nu_A(X)\in A\mbox{-inj}$ and $\nu_A(E^{-i}_X)\in\add(\nu_A(E))\subseteq\add(_AT)$ for all $i\geq 1$.


On the one hand, if $s\geq m\geq 1$, then $H^{-i}\big(\cpx{\Hom}_A(T, \nu_A(\cpx{E}_X))\big)=0$ for $0\leq i\leq m-1$. In this case, the following sequence
$$\nu_A(E^{-m}_X)\lra \nu_A(E^{-m+1}_X)\lra \cdots\lra \nu_A(E^{-1}_X)\lra \nu_A(X)\lra 0$$
can be regarded as a starting part of a minimal right $\add(_AT)$-approximation sequence of $\nu_A(X)$. Hence, up to isomorphism of complexes, we can write $\cpx{T}_X$ in the following form $$\cpx{T}_X:\quad \cdots\lra  T^{-m-1}_X\lra \nu_A(E^{-m}_X)\lra \nu_A(E^{-m+1}_X)\lra \cdots\lra \nu_A(E^{-1}_X)\lra \nu_A(X)\lra 0.$$
This leads to $\partial_T(X)\ge m$.

On the other hand, if $ \partial_T(X)\ge m\ge 1$, then $T^{-i}_X\in \add(\nu_A(E))\subseteq A\imodcat$ for all $1\leq i\leq m$.  Observe that, for any $I\in A\imodcat$, we have $D\Hom_A(\nu_A^{-}(I),-)\simeq \Hom_A(-,I)$ on $A\modcat$. Since the complex $\cpx{\Hom}_A(T,\cpx{T}_X)$ is exact, we apply $\nu_A^{-}$ to the complex $\cpx{T}_X$ and obtain the following sequence
$$ \nu_A^{-}(T^{-m}_X)\lra \nu_A^{-}(T^{-m+1}_X)\lra \cdots\lra \nu_A^{-}(T^{-1}_X)\lra X\lra 0$$
which induces an exact sequence of $B\opp$-modules:
$$
0\lra \Hom_A(X,T)\lra \Hom_A(\nu_A^{-}(T^{-1}_X),T)\lra \cdots\lra\Hom_A(\nu_A^{-}(T^{-m+1}_X),T)\lra \Hom_A(\nu_A^{-}(T^{-m}_X), T).
$$
Since $\nu_A^{-}(T^{-i}_X)\in \add(E)$ for each $1\leq i\leq m$, we see from Lemma \ref{proj-inj} that $\dm(\Hom_A(X,T)_B)\geq m$, that is, $s\ge m$. This finishes the proof of the first equality. $\square$

\medskip
Now, let us state a lower bound for the dominant dimension of $B$ in terms of $T$-gradients of $A$.

\begin{Koro}\label{equality}
$(1)$  $\dm(B)\geq \partial(A,T)\geq \partial_T(A)$.

$(2)$ $\partial_T(A)=\dm(T_B)=\nu_A(E)\emph{-}\dm(_AT)$.
\end{Koro}

{\it Proof.}
(1) Since $_AT$ is a tilting module, we see that $\Ext_A^j(T,T)=0$ for any $j\geq 1$. So the given projective resolution of $_AT$ induces a long exact sequence of $B^{\opp}$-modules:
$$
0\lra B_B \lra \Hom_A(P_0,T)\lra\Hom_A(P_1,T)\lra\cdots\lra \Hom_A(P_{n-1},T)\lra \Hom_A(P_n,T)\lra 0.
$$
Applying Corollary \ref{calculating dd}(1) to this sequence, we have
$$\dm(B_B)\geq \min\{\dm(\Hom_A(P_i,T)_B)+i\mid 0\leq i\leq n\}.$$
Thus $\dm(B)=\dm(B_B)\geq \min\{\partial_T(P_i)+i\mid 0\leq i\leq n\}=\partial(A,T).$ Since $\partial_T(P_i)\geq \partial_T(A)$ for all $0\leq i\leq n$, we clearly have $\partial(A,T)\geq \partial_T(A)$.

$(2)$ If we take $X={}_AA$ in Lemma \ref{approx}, then $\partial_T(A)=\dm(T_B)$.
Since $_AE$ is projective, there is an idempotent element $e=e^2\in A$ such that $\add(E)=\add(Ae)$.
It follows from Lemma \ref{approx} and Lemma \ref{idempotent} that  $\partial_T(A)=\nu_A(E)\emph{-}\dm(_AT)$.
$\square$

\subsection{Dominant dimensions for special classes of tilting modules \label{sect3.4}}

In this subsection we consider two special classes of tilting modules and estimate the dominant dimensions of their endomorphism algebras in
terms of their projective dimensions. Here considerations are focused on $n$-BB-tilting and canonical tilting modules.
The two classes of tilting modules have a common feature that all $P_i$ but one in their projective resolutions of $T$ belong to $\add(P)$. So we start with the following slightly general case.

\begin{Lem}\label{compare}
Suppose that $P_i\in\add(P)$ for all $0\leq i < n$. Then

$(1)$ $\dm(B)=\min\big\{\partial_T(P),\; n+\partial_T(P_n)\big\}.$

$(2)$ If $\dm(B)\geq n+1$, then $\dm(B)\leq n+\partial_T(A)$.
In this case, $\dm(B)=n+\partial_T(A)$ if and only if $\partial_T(P)\geq n+\partial_T(P_n)$.
\end{Lem}

{\it Proof.} As $_AT$ is a tilting module, we know that $\Ext_A^j(T,T)=0$ for any $j\geq 1$ and that the given projective resolution of $_AT$ induces a long exact sequence of $B^{\opp}$-modules:
$$
(\ast\ast)\quad 0\lra B_B \lra \Hom_A(P_0,T)\lra\cdots\lra \Hom_A(P_{n-1},T)\lra \Hom_A(P_n,T)\lra 0,
$$
with $\Hom_A(P_i,T)\in\add(\Hom_A(P,T))$ for all $0\leq i< n$ by the assumption $P_i\in\add(P)$.

Let $s:=\partial_T(P)$, $s':=\partial_T(P_n)$ and $t:=\dm(B)$. Then, since both $P$ and $P_n$ are projective, it follows from Lemma \ref{approx} that $s=\dm(\Hom_A(P,T)_B)$ and $s'= \dm(\Hom_A(P_n,T)_B)$. Consequently,
we have $\dm(\Hom_A(P_i,T))\geq s$ and $\dm(\Hom_A(P_0,T))=s$ since $\add(P)=\add(P_0)$ by assumption. Recall that
$$t=\dm(B^{\opp})=\dm(B_B)=\min\{\dm(Q_B)\mid Q\in\add(B_B)\}.$$
Since $P\in\add(_AT)$, we obtain $\Hom_A(P,T)\in\add(B_B)$. Thus $s\geq t$.

To show Lemma \ref{compare}, we shall apply Corollary \ref{calculating dd} to the sequence $(\ast\ast)$ by taking
$$\Lambda:=B^{\opp},\;m:=n,\;\, Y_{-1}:=B \;\, \mbox{and }\; \; Y_i:=\Hom_A(P_i,T)\; \mbox{for} \; 0\le i\le n.$$

Observe that $\dm(Y_i)\geq s\geq t$ for each $0\leq i\leq n-1$ and that $\dm(Y_0)=s$.
By Corollary \ref{calculating dd}(1), we first have  $t\geq \min\{s, n+s'\}$. By Corollary \ref{calculating dd}(2), we then obtain $s'\geq t-n$, and therefore $t\leq \min\{s, n+s'\}$. Thus $t=\min\{s, n+s'\}$. This finishes the proof of $(1)$.

$(2)$ Suppose $t\geq n+1$. 
By Corollary \ref{equality}(2) and $\add(\bigoplus_{i=0}^{n}P_i)=\add(P\oplus P_n)=\add(_AA)$, we see that
$\partial_T(A)= \dm(T_B)= \min\{\dm(\Hom_A(P,T)),\dm(\Hom_A(P_n,T))\}= \min\{s,s'\}$.
By $(1)$, we know that $t=\min\{s,n+s'\}\leq n+s'$. If $s\geq s'$, then $\partial_T(A)=s'$ and
$t\leq n+s'$. If $s<s'$, then $t=s=\partial_T(A)$. Thus we always have
$t\leq n+\partial_T(A)$, and the equality holds if and only if $s\geq n+s'$.  $\square$

\medskip
Lemma \ref{compare} can be applied to bound the dominant dimensions of the endomorphism algebras of a class of tilting modules by their projective dimensions. First of all, we mention the following technical result.

\begin{Koro}\label{small}
Suppose that $P_i\in\add(P)$ for $0\leq i<n$.

$(1)$ Let $f:E'\to P_n$ be a minimal right $\add(E)$-approximation of $P_n$, where $E'\in\add(E)$. If $\Hom_A(\Coker(f),T)\neq 0$, then $\dm(B)\leq n$.

$(2)$ Suppose that $_AE$ is injective. If $\omega \not\ncong E\,$ or $\,\nu_A(E)\not\ncong E$, then $\dm(B)\leq n$.
\end{Koro}

{\it Proof.} (1) If $\Hom_A(\Coker(f),T)\neq 0$, then $\partial_T(P_n)=0$ by Lemma \ref{approx}. Now, (1) follows from Lemma \ref{compare}(1).

(2) Suppose $\dm(B)\geq n+1$. By Lemma \ref{compare}(2), we have $\partial_T(A)\geq 1$. Since $\partial_T(A)=\nu_A(E)\emph{-}\dm(_AT)$ by Corollary \ref{equality}(2), there is an injection
$T\to I_0$ such that $I_0\in\add\big(\nu_A(E)\big)$. Since $\add(\omega)\subseteq\add(_AT)$, we have $\add(\omega)\subseteq \add\big(\nu_A(E)\big)$. By assumption, $_AE$ is injective, and therefore $E$ is projective-injective. Thus $\add(E)\subseteq \add(\omega)\subseteq \add\big(\nu_A(E)\big)$. However, both $\nu_A(E)$ and $E$ are basic and have the same number of indecomposable direct summands. This implies that  $\add(E)=\add(\omega)=\add\big(\nu_A(E)\big)$. Therefore $E\simeq \omega\simeq \nu_A(E)$. This proves (2). $\square$

\medskip
Now, we apply our results to $n$-BB-tilting modules which can be constructed, for instance, from Auslander-Reiten sequences (see \cite{hx2}).

Let $S$ be a simple, non-injective  $A$-module. Following \cite[Section 4]{hx2}, for an integer $n\ge 1$, we say that  $S$ defines an $n$-BB-tilting module if $\Ext_A^i(D(A),S)=0=\Ext_A^{i+1}(S,S)$ for $0\leq i\leq n-1$.
In this case,  we can associate a tilting module with $S$ in the following way:

Let $P(S)$ be the projective cover of $S$, and let $Q$ be the direct sum of all non-isomorphic indecomposable projective $A$-modules which are not isomorphic to $P(S)$. In \cite[Lemma 4.2]{hx2}, it is shown that the module $Q\oplus \tau^{-1}\Omega_A^{-n+1}(S)$ is a tilting $A$-module of projective dimension $n$, where $\tau^{-1}:=\rm{Tr}D$ is the Auslander-Reiten inverse translation. This tilting module is called the \emph{$n$-BB-tilting module} defined by $S$.
If, in addition, $S$ is projective, then this module will be called the \emph{$n$-APR-tilting module} defined by $S$.

\medskip
\begin{Koro}\label{BB-tilting}
Suppose that $_AT:=Q\oplus \tau^{-1}\Omega_A^{-n+1}(S)$ is an $n$-BB-tilting $A$-module.

$(1)$ If $\Hom_A(S,T)\neq 0$, then $\dm(\End_A(T))\leq n$.

$(2)$ Suppose that $S$ is projective. Then $\dm(\End_A(T))\leq n$. Moreover, if the injective envelope of $S$ is not projective, then
$\dm(A)\le 2n$.
\end{Koro}

{\it Proof.} (1) We first recall some properties of $n$-BB-tilting modules. Let
$$
0\lra S\lra I(S)\lra I_1\lra I_2\lra\cdots \lra I_{n-1}\lra I_n\lra I_{n+1}\lra\cdots
$$
be a minimal injective resolution of $S$, where $I(S)$ is the injective envelope of $S$.
Then $I_i\in\add(\nu_A(Q))$ for $1\leq i\leq n$ since $\Ext_A^i(S,S)=0$. Let $T':=\tau^{-1}\Omega_A^{-n+1}(S)$. Since the injective dimension of $S$ is at least $n$, we see that $0\neq T'$ does not contain projective direct summands. Applying $\nu_A^{-}$ to the above sequence, we obtain a minimal projective resolution of the module $T'$:
$$
0\lra P(S)\lra \nu_A^{-}(I_1)\lra \nu_A^{-}(I_2)\lra\cdots \lra \nu_A^{-}(I_{n-1})\lra \nu_A^{-}(I_n)\lra T'\lra 0
$$
such that $\nu_A^{-}(I_i)\in\add(Q)$ for $1\leq i\leq n$. This is due to $\Ext_A^j(D(A),S)=0$ for $0\leq j\leq n-1$.

For the module $T$, we recall that $E$ denotes a basic $A$-module such that $\add(_AE)=\{X\in\add(_AQ)\mid \nu_A(X)\in\add(_AT)\}.$ Further, let $f:E'\to P(S)$ be a minimal right $\add(E)$-approximation of $P(S)$ with $E'\in\add(E)$. Then the map $f$ cannot be surjective since $P(S)\notin\add(Q)$. This implies that the top of $\Coker(f)$ is equal to $S$. If
$\Hom_A(S,T)\neq 0$, then $\Hom_A(\Coker(f),T)\neq 0$, and therefore $\dm(\End_A(T))\leq n$ by Corollary \ref{small}(1).

(2) Note that if $S$ is projective, then $S=P(S)$ and $\Hom_A(S, Q)\neq 0$, because there exists an injection from $P(S)$ to $\nu_A^{-}(I_1)$. Since $\nu_A^{-}(I_1)\in\add(Q)\subseteq \add({_A}T)$, we have $\Hom_A(S,T)\neq 0$, and therefore $\dm(\End_A(T)) \le n$ by (1). If the injective envelope $I(S)$ of $S$ is not projective, then $\nu_A(S)= I(S)$ is not projective-injective. Consequently,
$S\not\in \add(\nu_A^{-}(\omega))\subseteq \add(Q)$. This means that $\nu_A^{-}(\omega)\in \add(E)$ and therefore $\omega\in \add(\nu_A(E))$. Thus $\dm(A) \le \dm(\End_A(T)) + n \le n+n=2n$ by  Theorem \ref{distance}(1).
$\square$

\medskip
Applying Corollary \ref{BB-tilting} to Auslander-Reiten sequences in module categories of Artin algebras, we obtain the following result.

\begin{Koro}
Let $0\to X\to M\to Y\to 0$ be an  Auslander-Reiten sequence in $A\modcat$ with $A$ an Artin algebra. If $X\ncong Y$, then $\dm\big(\End_A(X\oplus M\oplus Y)\big)\leq 2$.
\end{Koro}
{\it Proof.} Note that $X$ is not in $\add(M)$. Since $X\ncong Y$ by assumption, we have
$X\not\in\add(M\oplus Y)$. Let $V:=X\oplus M$, $U:=X\oplus M\oplus Y$ and $\Sigma:=\End_A(U)$.
Further, let $S^X$ be the top of the projective $\Sigma^{\opp}$-module $\Hom_A(X,U)$. By the remark (2) following \cite[Proposition 4.3]{hx2}, the $\Sigma^{\opp}$-module
$W:=\Hom_A(V,U)\oplus S^X$ is a $2$-BB-tilting module.  Let $\Delta:=\End_{\Sigma^{\opp}}(W)^{\opp}$. Then
$W$ is a $2$-APR-tilting $\Delta$-module defined by the simple projective $\Delta$-module $\Hom_{{\Sigma}^{\opp}}(S^X, W)$. Moreover, $\End_{\Delta}(W)\simeq \Sigma$ as algebras.
Now, it follows from Corollary \ref{BB-tilting}(2) that $\dm(\Sigma)\leq 2$. $\square$

\medskip
Next, we utilize our previous results to the so-called canonical tilting modules defined as follows:

Let $A$ be an algebra of $\dm(A)= n\ge 1$ and with a minimal injective resolution
$$0\lra {}_AA\lraf{d_0} E_0\lraf{d_1} E_1\lra \cdots \lraf{d_{n-1}} E_{n-1}\lra \cdots .$$
Let $T_i:=E_0\oplus \Coker(d_{i-1})$ and $B_i:=\End_A(T_i)$ for $1\le i<n+1$. It is not difficult to check that $T_i$ is a tilting $A$-module of projective dimension at most $i$ (for $i=1$, see also \cite[Proposition 5]{colby}). So these $T_i$ are called \emph{canonical tilting modules}.

Note that $A\PI \, = \add(\omega)=\add(E_0)$. For canonical tilting modules, we have the following proposition.

\begin{Prop} \label{invariant} With the above notation, we have the followings:

$(1)$ If $\add(E_0)\neq \add\big(\nu_A(E_0)\big)$, then $\dm(B_i)\leq i$ for all $1\le i<n+1$.

$(2)$ If $\add(E_0)=\add\big(\nu_A(E_0)\big)$, or equivalently, $\add(\soc(E_0))\simeq \add(\top(E_0))$, then $\dm(B_i)=n$ for all $1\le i<n+1$.

$(3)$ If $n$ is finite and $\nu_A(E_0)\in\add(T_n)$, then $\dm(B_n)=n$.
\end{Prop}

{\it Proof.}
Note that if the algebra $A$ is self-injective, then $E_0=A$ and $E_j=0$ for any $j\geq 1$. In this case, all the statements in Proposition \ref{invariant} are trivial. So, we now assume that $A$ is not self-injective. Since $\dm(A)=n\geq 1$, we always have $\add(\omega) = \add(E_0)$.

$(1)$ Working with a minimal injective resolution of $_AA$, we see that $\Omega^{-i}(A)$ contains no non-zero projective direct summands. Let $A=E_0'\oplus Q$, where $E_0'$ is injective and $Q$ does not contain any non-zero injective direct summands. Then, the module $T_i$ has a minimal projective resolution of the form
$$0\lra {}_AQ\lra E_0''\lra E_1\lra \cdots\lra E_{i-2}\lra E_{i-1}\oplus E_0\lra T_i\lra 0$$
where $E_0'\oplus E_0''\simeq E_0$. In particular, the module $E_{i-1}\oplus E_0$ is a projective cover of $T_i$. Let $W_i:=E(A,T_i)$ be the heart of $T_i$ (see Definition \ref{heart}). Then
$\add(W_i)=\{X\in\add(E_0)\mid \nu_A(X)\in\add(_AT_i)\}$ by definition.
Note that $W_i$ is injective and that $E_j\in\add(E_0)$ for $0\leq j\leq i-1$. By Corollary \ref{small}(2), if $\dm(B_i)\geq i+1$, then $\omega\simeq W_i\simeq \nu_A(W_i)$. Since $\add(\omega)=\add(E_0)$, we have
$\add(E_0)=\add\big(\nu_A(E_0)\big)$. This contradicts to the assumption in (1). Hence $\dm(B_i)\leq i$.

$(2)$ Suppose $\add(E_0)=\add\big(\nu_A(E_0)\big)$. Since $E_0\in\add(_AT_i)$, we  have $\add(W_i)=\add(E_0)=\add\big(\nu_A(W_i)\big)$. This implies $\partial_{T_i}(E_0)=\infty$. By Lemma \ref{compare}(2), we have $\dm(B_i)=i+\partial_{T_i}(Q)$. Note that
$$\partial_{T_i}(Q)=\min\big\{\partial_{T_i}(E_0),\partial_{T_i}(Q)\big\}=\partial_{T_i}(A)
$$
since $\add(E_0\oplus Q)=\add(_AA)$ and $\partial_{T_i}(E_0)=\infty$.
Thus $\dm(B_i)=i+\partial_{T_i}(A)$.
By Corollary \ref{equality}(2), we have $\partial_{T_i}(A)=\nu_A(W_i)$-$\dm(T_i)$. However, $\nu_A(W_i)\emph{-}\dm(_AT_i)=\dm(_AT_i)$ since $\add\big(\nu_A(W_i)\big)=A\PI$.
Hence, $$\dm(B_i)=i+\dm(_AT_i)=i+\dm(\Omega^{-i}_A(A))=n.$$

$(3)$ Suppose $\nu_A(E_0)\in\add(T_n)$. Then $\add(W_n)=\add(E_0)$ and therefore  $\partial_{T_n}(E_0)=\infty$. By Lemma \ref{compare}(2), we have $\dm(B_n)\geq n$. At the same time, it follows from $(1)$ and $(2)$ that $\dm(B_n)\leq n$. Thus $\dm(B_n)=n$.
$\square$

\medskip
In the following, we give an explicit description of the dominant dimension of $B_i$ for the case $\add(E_0)\neq \add\big(\nu_A(E_0)\big)$.

\begin{Koro} \label{Gen-Gorenstein}
Let $A$ be an algebra of finite dominant dimension $n\geq 1$ such that $\add(E_0)\neq \add\big(\nu_A(E_0)\big)$. Let $\mathcal{X}$, $\mathcal{Y}$ and $\mathcal{Z}$ be the complete sets of isomorphism classes $M$ of indecomposable, projective-injective $A$-modules such that $\nu_A(M)$ has projective dimension $0$, $n$ and $\geq n+1$, respectively. For each $M\in\mathcal{Y}\dot{\cup} \mathcal{Z}$, let
$$\cdots \cdots\lra Q^{-n-1}_M \lra Q^{-n}_M \lra Q^{-n+1}_M \lra \cdots\lra Q^{-1}_M\lra Q^{0}_M\lra \nu_A(M)\lra 0$$
be a minimal projective resolution of $\nu_A(M)$, and let $\delta_M$ denote the projective dimension of $\nu_A(M)$.

Suppose that $\delta_M<\infty$ and $Q_{M}^{-j}\in\add(E_0)$ for all $M\in\mathcal{Z}$ and
for all $n\leq j\leq \delta_M-1$.
Then
$$\dm(B_i)=\left\{\begin{array}{ll} \min\big\{\;i,\,\sigma(M,i)\mid M\in \mathcal{Y}\cup\mathcal{Z}\big\} & \mbox{if}\; 1\leq i\leq n-1,\\
\smallskip
\min\big\{n,\, \sigma(M,n)\mid M\in\mathcal{Z}\big\} & \mbox{if}\; i=n,\end{array} \right.$$
where
$$
\sigma(M,i):=\max_{0\leq\, j\leq\, \delta_M-i }\big\{j\mid Q_M^{-m}\in\add\big(\bigoplus_{X\in\mathcal{X}}\nu_A(X)\big)\;\;\mbox{for all}\;\; 0\leq m\leq j-1\big\}.
$$
In particular, if $\mathcal{Z}=\emptyset$, then $\dm(B_n)=n$ and $\dm(B_i)\leq \min\{i,n-i\}$ for $1\leq i\leq n-1$.
\end{Koro}

{\it Proof.} Recall that $A\PI\, =\add(E_0)$, due to $\dm(A)=n\geq 1$. Since $\add(E_0)\neq \add\big(\nu_A(E_0)\big)$, we see from Proposition \ref{invariant}(1) that $\dm(B_i)\leq i$ for $1\leq i\leq n$.  By Lemma \ref{compare}, we have $$\dm(B_i)=\min\{i,\,\partial_{T_i}(E_0)\}.$$  Next, we shall describe $\partial_{T_i}(E_0)$ explicitly.

Since $\dm(A)=n=\dm(A^{\opp})$, each indecomposable injective, non-projective $A$-module $V$ has
projective dimension at least $n$ with a minimal projective resolution of the form
$$
\cdots \lra Q_V^{-n} \lra Q_V^{-n+1} \lra \cdots\lra Q_V^{-1}\lra Q_V^{0}\lra V\lra 0
$$
such that $Q_V^{-j}\in\add(E_0)$ for $0\leq j\leq n-1$. In particular, if $\pd(V)=n$, then
$\Omega_A^n(V)\in\add(_AA)$, and therefore $V\simeq \Omega_A^{-n}(\Omega_A^n(V))\in\add\big(\Omega_A^{-n}(A)\big)=\add\big(\Coker(d_{n-1})\big)$. This implies that if $M\in\mathcal{Y}$, then we automatically have $Q_{M}^{-j}\in\add(E_0)$ for all $0\leq j\leq \delta_M-1=n-1$, and therefore $\nu_A(M)\in \add(\Coker(d_{n-1}))$.

Let $W:=\bigoplus_{X\in\mathcal{X}}X$. Then $W$ is a basic $A$-module such that
$$
\add(W)=\{X\in\add(E_0)\mid \nu_A(X)\in\add(E_0)\}\quad \mbox{and}\quad
\add(W\oplus \bigoplus_{M\in\mathcal{Y}\dot\cup\mathcal{Z}}M)=\add(E_0).
$$
Moreover, the hearts of the canonical tilting $A$-modules $T_i$ are given by
$$
E(A, T_s) = W \;\;\mbox{for}
\;\;1\leq s\leq n-1\;\;\mbox{and}\;\;
E(A,T_n)= W\oplus \bigoplus_{M\in\mathcal{Y}}M.
$$
This implies that
$\partial_{T_s}(W)=\infty=\partial_{T_n}\big(W\oplus \bigoplus_{M\in\mathcal{Y}}M\big)$ with $1\le s\le n-1$.
Thus $$\partial_{T_s}(E_0)=\partial_{T_s}(\bigoplus_{M\in\mathcal{Y}\dot\cup\mathcal{Z}}M)
=\min_{M\in\mathcal{Y}\dot\cup\mathcal{Z}}\{\partial_{T_s}(M)\}
\quad \mbox{and}\quad
\partial_{T_n}(E_0)=\partial_{T_n}(\bigoplus_{M\in\mathcal{Z}}M)
=\min_{M\in\mathcal{Z}}\{\partial_{T_n}(M)\}.$$
To show Corollary \ref{Gen-Gorenstein}, we only need to prove the following two statements:

$(a)$ $ \partial_{T_s}(M)=\sigma(M,s)$ for $1\le s\le n-1$ and any $M\in\mathcal{Y}\dot\cup\mathcal{Z};$

$(b)$ $ \partial_{T_n}(M)=\sigma(M,n)$ for any $M\in\mathcal{Z}$.

\smallskip
\noindent Indeed, for any $M\in\mathcal{Y}\dot\cup\mathcal{Z}$, we have $\pd(\nu_A(M))=\delta_M<\infty$ and $Q_{M}^{-j}\in\add(E_0)$ for $0\leq j\leq \delta_M-1$. This implies that $Q^{-\delta_M}_M\in\add(_AA)$ and  $\Omega_A^{-i}(Q^{-\delta_M}_M)\in\add(\Omega_A^{-i}(A))\subseteq\add(T_i)$ for all $1\leq i\leq n$. So the following sequence
$$0\lra \Omega_A^{-i}(Q^{-\delta_M}_M)\lra  Q^{-(\delta_M-i-1)}_M \lra \cdots\lra Q^{-1}_M\lra Q^{0}_M\lra \nu_A(M)\lra 0,$$
is a minimal right $\add(T_i)$-approximation sequence of $\nu_A(M)$. Note that $\Omega_A^{-i}(Q^{-\delta_M}_M)$ is injective if and only if $i=n$ and $M\in\mathcal{Y}$. Thus  $(a)$ and $(b)$ hold. This finishes the proof of Corollary \ref{Gen-Gorenstein}. $\square$

\medskip
We observe that all Morita algebras satisfy the condition in Proposition \ref{invariant}(2).

\begin{Koro} Let $A$ be a Morita algebra with $E_0$ the injective envelope of $_AA$.
Then $\End_A(E_0\oplus
\Omega^{-j}(A))$ is again a Morita algebra and has dominant dimension equal to
$\dm(A)$ for all $1\le j<
\dm(A)+1.$ \label{cor23}
\end{Koro}

{\it Proof.} Since $A$ is a Morita algebra, it follows from Lemma \ref{Morita} that $\add(E_0)=A\PI\, =\mathscr{E}(A)$. By Lemma \ref{stp}(1),  we further have $\add(E_0)=
\add\big(\nu_{A}(E_0)\big)$. Thus $\dm\big(\End_A(E_0\oplus
\Omega^{-j}(A))\big)=\dm(A)\ge 2$ by Proposition \ref{invariant}(2). It follows from $\add(E_0)=\add\big(\nu_A(E_0)\big)$ that $\mathscr{E}(A)=\add(E)$ where $E$ is the herat of $T_j:=E_0\oplus
\Omega^{-j}(A)$. Thus, by Corollary \ref{B-M}, the algebra $\End_A(T_j)$ is a Morita algebra. $\square$

\medskip
Finally, we use an example in \cite{FH} to illustrate Corollary \ref{Gen-Gorenstein} and show that the equality $``\dm(B_i)=\min\{i,n-i\}"$ really occurs for $1\leq i\leq n-1$.

Let $n\geq 3$ and let $A$ be the quotient of the path algebra over a field $k$ of the following quiver
$$2n+1\lraf{\beta_{2n}} 2n\lraf{\beta_{2n-1}}\cdots\lra n+2\lraf{\beta_{n+1}} n+1\lraf{\beta_n} n\lraf{\beta_{n-1}} n-1\lra\cdots\lraf{\beta_3} 3\lraf{\beta_2} 2\lraf{\beta_1} 1$$
with relations $\beta_{i+1}\beta_i=0$ for $1\leq i\leq 2n-1$ except $i=n$. Further, let $P(i)$, $I(i)$ and $S(i)$ be the indecomposable projective, injective and simple $A$-modules corresponding
to the vertex $i$ for $1\leq i\leq 2n+1$, respectively.

Let $W:=\bigoplus_{i=2}^{n}(P(i)\oplus P(n+i))$. Then
$$
A\PI\, =\add(W\oplus P(2n+1))\;\;\mbox{and}\;\;\add(W)=\{X\in A\PI\, \mid \nu_A(X)\in A\PI\,  \}.
$$
In fact, we have $\nu_A(P(j))\simeq P(j+1)$ for $1\leq j\leq n-1$ and $n+2\leq j\leq 2n$, $\nu_A(P(n))=P(n+2)$
and $\nu_A(P(2n+1))\simeq I(2n+1)$. 
Note that $_AA$ has a minimal injective resolution of the form:
$$0\lra {}_AA\lra E_0\lra E_1\lra \cdots \lra E_{n-1}\lra E_n\lra 0$$
where $E_0:=(W\oplus P(2n+1))\oplus P(2)\oplus P(n+2)$,\, $E_i:=P(i+2)\oplus P(n+i+2)$ for $1\leq i\leq n-2$,\, $E_{n-1}:=P(n+2)\oplus P(2n)$ and $E_n:=I(n+1)\oplus I(2n+1).$
Thus $\dm(A)=n$.

Clearly, $\add(E_0)\ncong \add(\nu_A(E_0))$, $\mathcal{X}=\{P(j), P(n+j)\mid 2\leq j\leq n \}$, $\mathcal{Y}=\{P(2n+1)\}$ and $\mathcal{Z}=\emptyset$ (see Corollary \ref{Gen-Gorenstein}). Since $\nu_A(P(2n+1))\simeq I(2n+1)$ has a minimal projective resolution
$$
0\lra P(n+1)\lra P(n+2)\lra P(n+3)\lra \cdots\lra P(2n-1)\lra P(2n)\lra P(2n+1)\lra I(2n+1)\lra 0
$$
and since $\nu_A(\mathcal{X})=\{P(s), P(t) \mid 3\leq s\leq n,\; n+2\leq t\leq 2n+1\}$, we have $\sigma({P(2n+1)},i)=n-i$ for $1\leq i\leq n-1$.

Now, let $T_i:=E_0\oplus\Omega_A^{-i}(A)$ and $B_i:=\End_A(T_i)$ for $1\leq i\leq n$, where $\Omega_A^{-i}(A)=S(i+1)\oplus S(n+i+1)\;\;\mbox{for}\;\;1\leq i\leq n-1\;\;\mbox{and}\;\; \Omega_A^{-n}(A)=I(n+1)\oplus I(2n+1).$ It follows from Corollary \ref{Gen-Gorenstein}
that $\dm(B_i)=\min\{i,n-i\}$ for $1\leq i\leq n-1$. Clearly, $\dm(B_n)=n$ since $\mathcal{Z}=\emptyset$.

This example also shows  that the condition $`` \nu_A(E)\in \add(\omega)$" does not imply $``\omega\in \add(\nu_A(E))"$, where $\omega$ is an additive generator for $A\PI $. Indeed, for
a fixed $T_i$ with $i\ne n$, we can verify that $\nu_A(E)\in \add(\omega)$ and $\omega\not\in \add(\nu_A(E))$. The last two conditions are equivalent to that $``\omega'\in \add(\nu_{B_i^{\opp}}(E'))$" and $``\, \nu_{B_i^{\opp}}(E')\not\in \add(\omega')"$, where $\omega'$ stands for an additive generator for $B_i^{\opp}\PI$ (see the proof of Lemma \ref{dm-formula}).

\section{Are generalized symmetric algebras invariant under derived equivalences? \label{sect4}}
As is known, many important algebras (for example, Schur algebras and
$q$-Schur algebras) are Morita equivalent to
algebras of the form $\End_A(A\oplus X)$ with $A$ a symmetric
algebra (that is, ${}_AA_A\cong {}_AD(A)_A)$ and $X$ an $A$-module. Algebras of this form were called \emph{generalized symmetric algebras} in \cite{FK}. Clearly, generalized symmetric algebras are Morita algebras. It is known that symmetric algebras (or more generally, self-injective algebras) are closed under derived equivalences (see \cite{Rickard}). So, when Ming Fang (with his coauthors) studies generalized symmetric algebras and certain
quasi-hereditary covers of some Hidden Hecke algebras
as well as the dominant dimensions of blocks of
$q$-Schur algebras, he asks naturally the following
question: Are  generalized symmetric algebras closed under
taking derived equivalences? More precisely,

\medskip
{\bf Question.} Let  $\Lambda$ and $\Gamma$ be finite-dimensional
$k$-algebras over a field $k$ such that they are derived equivalent.
Suppose that $\Lambda$ is of the form $\End_A(A\oplus X)$ with $A$ a
symmetric $k$-algebra and $X$ an $A$-module. Is there a symmetric
$k$-algebra $B$ and a $B$-module $Y$ such that $\Gamma$ is
isomorphic to $\End_B(B\oplus Y)$?

\medskip
As mentioned, for a symmetric algebra $A$, if $X=0$, then the
above question gets positive answer (see \cite{Rickard}). In this
section, however, we shall apply our results in the previous sections to give a
negative answer to the above question for $X\ne 0$.


\subsection{A negative answer}
Our concrete counterexample to the
above question is actually another consequence of Theorem \ref{main result}.

\begin{Koro}\label{cor2} Suppose that $k$ is a field with an element that is not a root of unity. Then there exist two finite-dimensional $k$-algebras $\Lambda$ and
$\Gamma$ satisfying the following conditions:

$(1)$ $\Lambda$ and $\Gamma$ are derived equivalent with the same
Cartan matrices. In particular, $\dk(\Lambda)=\dk(\Gamma)$.

$(2)$ There is a symmetric algebra $A$ and a generator $M$ over $A$
such that $\Lambda\simeq \End_A(M)$ as algebras.

$(3)$ $\dm(\Lambda)=2$ and $\dm(\Gamma)=1$.

\end{Koro}

In Corollary \ref{cor2}, since $\dm(\Gamma)=1$, the algebra $\Gamma$
cannot be Morita equivalent to a generalized symmetric algebra.

Our discussion in the sequel is partially based on some results in \cite{CPX}. So, we first
recall necessary ingredients from \cite{CPX}.

From now on, we fix a non-zero element $q$ in a fixed field $k$, and
suppose that $q$ is not a root of unity.

Let $A$ be the Liu-Schulz $k$-algebra in \cite{Ls}, that is, $A$ is an
$8$-dimensional unitary $k$-algebra with the generators:
$x_0, x_1, x_2$, and the relations: $ x_i^2= 0  \mbox{ and }
x_{i+1}x_i+q x_ix_{i+1}= 0 \mbox{ for } i= 0, 1, 2,$ where the
subscripts are modulo $3$. Note that $A$ is a local, symmetric $k$-algebra.

For $j\in \mathbb{Z}\,$, we set
$$u_j: =x_2+q^jx_1,\quad  I_j := Au_j \quad \mbox{and}\quad \Lambda_j:=
\End_A(A \sz I_0 \sz I_j ).$$

For the algebra $A$ and
the modules $I_j$, we cite the following properties from \cite[Section 6]{CPX}.

\begin{Lem}\label{ls}
Let $i$ and $j$ be integers. Then the following statements are true:

$(1)$  As $k$-algebras, $A/\rad(A)\simeq k\simeq
\End_A(I_i)/\rad\big(\End_A(I_i)\big)$.

$(2)$ If $i\neq j$, then $I_i$ and $I_j$ are non-isomorphic as
$A$-modules.

$(3)$ There exists a short exact sequence $0\ra I_{i+1}\ra A \ra
I_i \ra 0$ of $A$-modules.

$(4)$ $\dk{I_i}=4=\dk{\Hom_A(I_i, A)}$ and $\dk{\Hom_A(I_j, I_i)}=
\left\{\begin{array}{ll} 3 & \mbox{if } j=i \mbox{ or } i-2,\\
2 & \mbox{otherwise}.\end{array} \right.$

$(5)$ $\dk{\Ext_A^1(I_j,I_i)}=
\left\{\begin{array}{ll} 1 & \mbox{if } j\leq i\leq j+3,\\
0 & \mbox{otherwise}.\end{array} \right.$
\end{Lem}

\medskip
Note that, by \cite[Proposition 6.9]{CPX}, all the algebras
$\Lambda_j$ for $j\geq 3$ are derived equivalent, but not stably
equivalent of Morita type. Now, let us look at their Cartan matrices. By Lemma \ref{ls}(1) and (4),
one can easily calculate the Cartan matrix $C_{\Lambda_2}$ of
$\Lambda_2$ and the Cartan matrix $C_{\Lambda_3}$ of $\Lambda_3$:
$$
C_{\Lambda_2}= \left(
   \begin{array}{ccc}
   8  & 4 & 4\\
   4  & 3 & 3 \\
   4  & 2 & 3 \\
   \end{array}
 \right)\quad \mbox{and}\quad
C_{\Lambda_3}= \left(
   \begin{array}{ccc}
   8  & 4 & 4\\
   4  & 3 & 2 \\
   4  & 2 & 3 \\
   \end{array}
 \right).
 $$
Clearly, the former is not a symmetric matrix, but the latter is.
Since the Cartan matrices of two derived equivalent algebras are
congruent over $\mathbb Z$, derived equivalences preserve the
symmetry of Cartan matrices. Thus $\Lambda_2$ is not
derived equivalent to $\Lambda_3$.

In the following, however, we shall show that $\Lambda_2$ has a subalgebra
$\Gamma$ which is derived equivalent to $\Lambda_3$ such that
$C_\Gamma=C_{\Lambda_3}$. This may illustrate an intrinsic connection between $\Lambda_2$ and $\Lambda_3$

\medskip
{\bf Proof of Corollary \ref{cor2}.} We
define $\Lambda:=\Lambda_3=\End_A(A\oplus I_0\oplus I_3)$ and
$$\Gamma:=\left(
 \begin{array}{ccc}
\End_A(A)      & \Hom_A(A, I_0) &   \Hom_A(A, I_2) \\
\Hom_A(I_0, A) & \End_A(I_0)    & \mathscr{P}(I_0, I_2)\\
\Hom_A(I_2, A)   & \Hom_A(I_2, I_0) &  \End_A(I_2)  \\
   \end{array}
 \right).$$
Clearly, $\Gamma$ is a subalgebra of
$\Lambda_2:=\End_A(A\oplus I_0\oplus I_2)$.

First, we prove that $\Lambda$ and $\Gamma$ have the same
Cartan matrix.

In fact, it is sufficient to check $\dk{\mathscr{P}(I_0, I_2)}=2$.
By Lemma \ref{ls}(3), there is an exact sequence $0\to
I_3\to A \to I_2 \to 0$ of $A$-modules. Applying $\Hom_A(I_0, -)$ to
this sequence, we get another exact sequence of $k$-modules:
$$0\lra \Hom_A(I_0, I_3)\lra \Hom_A(I_0, A) \lra  \mathscr{P}(I_0, I_2) \lra 0.$$
Note that $\dk{\Hom_A(I_0, I_3)}=2$ and $\dk{\Hom_A(I_0, A)}=4$ by
Lemma \ref{ls}(4). Thus $\dk{\mathscr{P}(I_0, I_2)}=2$. This
implies that $\Gamma$ and $\Lambda$ have the same Cartan matrices. Clearly, $\dk{\Lambda}=\dk{\Gamma}=34.$

Second, we take $N:=I_0$ and $Y:=I_2$ in Corollary \ref{cor1}. It
follows from Lemma \ref{ls}(3) that $I_3=\Omega_A(I_2)\,$, and form
Lemma \ref{ls}(5) that $\Ext^1_A(I_2, I_0)=0$ and $\Ext^1_A(I_0,
I_3)\neq 0$. Now, by Corollary \ref{cor1}, we know that $\Lambda$
and $\Gamma$ are derived equivalent with $\dm(\Gamma)=1$.

Since $\Ext_A^1(I_0, I_0)\neq 0$ by Lemma \ref{ls}(5), we have
$\dm(\Lambda)=2$ by \cite[Lemma 3]{Muller} (or the remarks after the definition of dominant dimension).

Consequently, $\Lambda$ and $\Gamma$ satisfy all the properties mentioned in Corollary \ref{cor2}. $\square$

\medskip
Observe that the algebras $\Lambda$ and $\Gamma$ defined in the
proof of Corollary \ref{cor2} can be described by matrices.

Let $V$ be a $k$-vector space. For $y_i\in V$ with $1\leq i\leq n\in
\mathbb{N}$, we denote by $<y_1,\ldots,y_n>$ the $k$-subspace
of $V$ spanned by all $y_i$.

Define $J_j:=u_j \,A$ for $j\in\mathbb{Z}$ and
$${\small C:=<1,x_1,x_2,x_0x_1,x_1x_2,x_2x_0,x_0x_1x_2>,\;\,
T:=<x_1,x_2,x_0x_1,x_1x_2,x_2x_0,x_0x_1x_2>, \;S:=T\oplus <x_0>}$$

It has been shown in \cite[Proposition 6.9]{CPX} that
$$
 \Lambda:=\Lambda_3\simeq \left(
   \begin{array}{ccc}
     A   & A/I_1 & A/I_4 \\
     J_0 & C/I_1 & T/I_4 \\
     J_3 & T/I_1 & C/I_4 \\
   \end{array}
\right)\quad \mbox{and}\quad \Lambda_2\simeq \left(
   \begin{array}{ccc}
     A   & A/I_1 & A/I_{3} \\
     J_0 & C/I_1 & S/I_{3} \\
     J_2 & T/I_1 & C/I_{3} \\
   \end{array}
\right)
$$
as algebras. By the proof of Corollary \ref{cor2}, the algebra
$\Gamma$ is a subalgebra of $\Lambda_2$, which, with help of
\cite[Lemma 6.3(4), Lemma 6.5]{CPX}, can be described as the following matrix algebra:
$$\Gamma\simeq
\left(
   \begin{array}{ccc}
     A   & A/I_1 & A/I_{3} \\
     J_0 & C/I_1 & T/I_{3} \\
     J_2 & T/I_1 & C/I_{3} \\
   \end{array}
\right).$$

Recall that, in the last part of \cite[Section 6]{CPX}, we have
described the algebra $\Lambda_m$ for any $m\geq 3$ (up to
isomorphism) by a fixed quiver $Q$ with relations $\rho_m$. Note
that the definition of $\rho_m$ depends on $m$ and makes sense for $m=2$. Thus we can show that $\Gamma$ is actually
isomorphic to the algebra $kQ/<\rho_2>$,
where $<\rho_2>$ is the ideal of the path algebra $kQ$
generated by the set of relations $\rho_2$.

\subsection{Further questions}
Finally, we mention the following conjectures suggested by the results in this paper.

\smallskip
(1) If two algebras $A$ and $B$ are derived equivalent (not necessarily given by a tilting module), then $\dm(A)=\infty$ if and only if $\dm(B)=\infty$. Equivalently, if $A$ and $B$ are derived equivalent, then $\dm(A)<\infty$ if and only if $\dm(B)<\infty$.

(2) If $A$ is an algebra of dominant dimension $n\ge 1$ with a minimal injective resolution
$ 0\ra {}_AA\ra E_0\ra \cdots \ra E_{n-1}\ra \cdots$, then $\dm\big(\End_A(E_0\oplus \Omega^{-n}(A))\big)=n$.

(3) Suppose that two algebras $A$ and $B$ are derived equivalent. If $A$ is a Morita algebra and $\dm(B)\ge 2$, then $B$ is a Morita algebra.

\bigskip
{\bf Acknowledgement.} The authors thank Ming Fang for conveying the question to them and
explaining backgrounds of the question. The research work is partially supported by BNSF and NNSF.

\bigskip
First version: March 25, 2012. Revised: July 6, 2014. Final version: January 31, 2015.
\end{document}